\documentclass [leqno, spanish, 10pt]{amsart}

\usepackage[T1]{fontenc}
\usepackage[utf8]{inputenc}

\usepackage[all]{xy}
\usepackage{amssymb, amsmath, amsfonts, amsthm, amsbsy, amscd, amsxtra, stmaryrd}
\usepackage[bbgreekl]{mathbbol} 
\usepackage[vcentermath]{youngtab} 
\usepackage[mathscr]{euscript}
\usepackage{upgreek}
\usepackage{mathtools}
\usepackage{tikz}
\usepackage{makecell} 
\usepackage{url}
\usepackage{hyperref}
\hypersetup{
	colorlinks=false, linktocpage=true, hidelinks = true
}

\usepackage{etoolbox}
\apptocmd{\sloppy}{\hbadness 10000\relax}{}{}

\usepackage{amsthm}
\newtheorem{theorem}{Theorem}[section]
\newtheorem{corollary}[theorem]{Corollary}
\newtheorem{lemma}[theorem]{Lemma}
\newtheorem{definition}	[theorem]{Definition}

\theoremstyle{definition}

\usepackage{thmtools}
\usepackage{lipsum}

\declaretheoremstyle[
  headfont=\normalfont\bfseries,
  sharenumber = theorem,
  bodyfont=\normalfont,
  qed={$\boxtimes$},
]{examplestyle2}

\declaretheorem[
  style=examplestyle2,
  title=Example,
  refname={example,examples},
  Refname={Example,Examples}
]{example}

\declaretheorem[
  style=examplestyle2,
  title=Remark,
  refname={remark, remarks},
  Refname={Remark, Remarks}
]{remark}

\numberwithin{equation}{section}

\newcommand{\imat}{I}
\newcommand{\kform}{\kappa}
\newcommand{\miya}{\tau}
\newcommand{\fusion}{\mathscr{F}}
\newcommand{\axes}{\mathscr{A}}
\DeclareMathAlphabet{\mathantt}{OT1}{antt}{li}{it}
\DeclareMathAlphabet{\mathpzc}{OT1}{pzc}{m}{it}
\renewcommand{\P}{\mathscr{P}}

\newcommand{\card}{\#}
\newcommand{\F}{\mathscr{F}}

\newcommand{\matsuo}{\mathbb{M}}
\newcommand{\ualg}{\hat{\talg}}
\newcommand{\talg}{\mathbb{S}}

\newcommand{\join}{\wedge}

\newcommand{\Id}{\operatorname{Id}}

\newcommand{\C}{\mathscr{C}} 
\newcommand{\ones}{\mathbb{1}}

\newcommand{\setn}{\bar{\mathsf{n}}}
\newcommand{\setm}{\bar{\mathsf{m}}}

\newcommand{\ealg}{\mathbb{E}}

\newcommand{\ideal}{\mathbb{I}}

\newcommand{\fie}{\mathbb{k}}

\newcommand{\imt}{\iota}

\newcommand{\alg}{\mathbb{A}}

\newcommand{\balg}{\mathbb{B}}

\newcommand{\om}{\omega}
\newcommand{\mlt}{\circ}

\newcommand{\hmlt}{\hat{\circ}}

\newcommand{\dum}{\,\cdot\,\,}

\newcommand{\la}{\lambda}
\newcommand{\ep}{\epsilon}
\newcommand{\fiet}{\fie^{\times}}

\newcommand{\eno}{\operatorname{End}}
\newcommand{\si}{\sigma}

\newcommand{\integer}{\mathbb{Z}}
\newcommand{\ztwo}{\integer/2\integer}

\newcommand{\B}{\mathscr{B}}
\newcommand{\lb}{\langle}
\newcommand{\ra}{\rangle}

\newcommand{\ste}{\mathbb{V}}

\newcommand{\al}{\alpha}
\newcommand{\be}{\beta}
\newcommand{\ga}{\gamma}

\DeclareMathOperator{\Aut}{Aut}

\newcommand{\rea}{\mathbb R}

\newcommand{\tr}{\operatorname{\mathsf{tr}}}

\let\oldtocsection=\tocsection
\let\oldtocsubsection=\tocsubsection
\renewcommand{\tocsection}[2]{\hspace{0em}\oldtocsection{#1}{#2}}
\renewcommand{\tocsubsection}[2]{\hspace{1em}\oldtocsubsection{#1}{#2}} 

\begin{document}
\title[Commutative algebras associated with Steiner triple systems]{Killing metrized commutative nonassociative algebras associated with Steiner triple systems}
\author{Daniel J.~F. Fox} 
\address{Departamento de Matemática Aplicada\\ Escuela Técnica Superior de Arquitectura\\ Universidad Politécnica de Madrid\\Av. Juan de Herrera 4 \\ 28040 Madrid España}
\email{daniel.fox@upm.es}

\begin{abstract}
With each Steiner triple system there is associated a one-parameter family of commutative, nonassociative, nonunital algebras that are by construction exact, meaning that the trace of every multiplication operator vanishes, and these algebras are shown to be Killing metrized, meaning the Killing type trace-form is nondegenerate and invariant (Frobenius), and simple, except for certain parameter values.  The definition of these algebras resembles that of the Matsuo algebra of the Steiner triple system, but they are different. 
For a Hall triple system, the associated algebra is a primitive axial algebra for a $\ztwo$-graded fusion law. 
\end{abstract}

\maketitle

\section{Introduction}
Let $(\alg, \mlt)$ be a finite-dimensional commutative algebra over a field $\fie$ of characteristic zero. The multiplication $\mlt$ is not supposed either associative or unital and in the examples considered here it will be neither. The class of commutative $\fie$-algebras is too large to admit an interesting theory, so it is natural to look for conditions that determine tractable classes of such algebras. Let $L:\alg \to \eno(\alg)$ be the multiplication operator defined by $L(x) = x\mlt y$. If $\Phi$ is an isomorphism of commutative algebras, then $L(\Phi(x)) = \Phi\circ L(x)\circ \Phi^{-1}$, so any polynomial on $\eno(\alg)$ invariant with respect to the automorphism group of $\Aut(\alg, \mlt)$ of $(\alg, \mlt)$ determines an isomorphism invariant of commutative algebras. In particular $\tr L(x_{1})\dots L(k_{k})$ is an invariant for all $k \geq 1$. See \cite{Popov} for general discussion of such trace forms from the perspective of invariant theory. A commutative algebra $(\alg, \mlt)$ is \emph{exact} if $\tr L(x) = 0$ for all $x \in \alg$ (note that an exact algebra is not unital). (Exact algebras are called \emph{harmonic} in \cite{Nadirashvili-Tkachev-Vladuts}.) Define the \emph{Killing} type trace-form $\kform$ of $(\alg, \mlt)$ by $\kform(x, y) = \tr L(x)L(y)$. A commutative algebra $(\alg, \mlt)$ is \emph{Killing metrized} if its Killing form $\kform$ is nondegenerate and invariant. In general it is interesting that a commutative algebra admit an invariant bilinear form (often said to be a \emph{Frobenius} form). It is a stronger condition that it be metrized by a specific trace-form, as such a form is determined by the underlying algebraic structure and so is not additional data. For an exact commutative algebra the only possibility is the Killing type form, for any nontrivial bilinear form constructed from traces of the multiplication endomorphisms and the multiplication itself must be a multiple of the Killing form. The class of Killing metrized exact commutative algebras is perhaps the simplest interesting class defined by conditions on trace-forms, and it contains the tensor products of pairs of semisimple Lie algebras, the deunitalizations of simple Euclidean Jordan algebras, and many more examples. For this and further motivation for considering Killing metrizability see \cite{Fox-simplicial, Fox-cubicpoly}.

Here there is associated with each Steiner triple system $\B$ on $\setn = \{1, \dots, n\}$ a one-parameter family of $(n-1)$-dimensional exact commutative nonassociative $\fie$-algebras $(\talg_{\be}(\fie, \B), \mlt)$. The algebras $(\talg_{\be}(\fie, \B), \hmlt)$ are quotients by a one-dimensional ideal of members of a $1$-parameter subfamily of the $3$-parameter family of algebras called \emph{triple algebras} in \cite[Section $5.3$]{Hall-transpositionalgebras}, which should be consulted for further context and references. Their description in terms of generators is given in \eqref{talgrelations}, following their Definition \ref{talgdefined}.

Theorem \ref{stsalgebratheorem} shows that $(\talg_{\be}(\fie, \B), \mlt)$ is Killing metrized for all $\be \in \fie$ such that $(n-3)\be^{2} + 1 \neq 0$. Moreover, when $\fie$ is a Euclidean field (a formally real field for which every nonzero element is a square or minus a square), $\kform$ is positive definite for all $\be \in \fie$. Theorem \ref{stsalgebrasimpletheorem} shows that $(\talg_{\be}(\fie, \B), \mlt)$ is simple for $\be \notin \{1,  -\tfrac{n-1}{n-3}\}$. Example \ref{notsimpleexample} shows that $(\talg_{1}(\fie, \B), \mlt)$ need not be simple; when $\B$ is the affine plane of order $3$ the $8$-dimensional algebra $(\talg_{1}(\fie, \B), \mlt)$ decomposes as a direct sum of $4$ mutually orthogonal, mutually isomorphic, $2$-dimensional subalgebras.

Section \ref{hallsection} describes stronger results that are obtained when $\B$ is a Hall triple system. Theorem \ref{htstheorem} shows that for $\be \in \fie \setminus\{-\tfrac{n-1}{n-3}, 0, 1\}$ the algebra $(\talg_{\be}(\fie, \B), \mlt)$ is a primitive axial algebra for a $\ztwo$-graded fusion law in the sense of \cite{DeMedts-Peacock-Shpectorov-VanCouwenberghe, Hall-Rehren-Shpectorov, Khasraw-McInroy-Shpectorov}, and that the Miyamoto group of this axial algebra structure acts as automorphisms of $\B$ and $(\talg_{\be}(\fie, \B), \mlt)$. (The relevant definitions are recalled in Section \ref{hallsection}.) When $\be = \tfrac{1}{n-1}$, the algebra $(\talg_{\be}(\fie, \B), \mlt)$ is a primitive axial algebra for a Jordan fusion law; in this case it is a quotient of a Matsuo algebra by a one-dimensional ideal. From the point of view of axial algebras, it seems interesting that the local affine condition underlying the definition of Hall triple system determines an axial structure. 

Section \ref{idempotentsection} records some partial results regarding the structure of idempotents in these algebras.

\section{Preliminaries and definition of the algebras}
Before defining $(\talg_{\be}(\fie, \B), \hmlt)$ and stating the main results, there is fixed the terminology regarding the combinatorial objects of interest. For background see \cite{Aschbacher, Colbourn-Dinitz, Dembowski}.

An \emph{incidence structure} $(\P, \B, \F)$ comprises sets $\P$ and $\B$ and an \emph{incidence relation} $\F \subset \P \times \B$. Elements $p \in \P$ and $B \in \B$ are \emph{incident} if $(p, B) \in \F$. The letters $\P$, $\B$, and $\F$ are used to suggest the words \emph{points}, \emph{blocks}, and \emph{flags}. Blocks are often called \emph{lines}. The \emph{dual} incidence structure $(\B, \P, \F^{\ast})$ is defined by the dual incidence relation $\F^{\ast} = \{(B, p): (p, B) \in \F\}$ of $\B \times \P$. A subspace of an incidence structure $(\P, \B, \F)$ is an incidence structure $(\bar{\P}, \bar{\B}, \bar{\F})$ such that $\bar{\P} \subset \P$, $\bar{\B} \subset \B$, and $\bar{\F} \subset \F$. Let $\pi_{\P}:\P \times \B \to \P$ and $\pi_{\B}:\P \times \B \to \B$ be the projections onto the factors.  A morphism of incidence structures $(\P_{i}, \B_{i}, \F_{i})$, $i \in \{1, 2\}$ is a map $\Phi:\F_{1} \to \F_{2}$ such that $\pi_{\P_{2}}\circ \Phi(\pi_{\P_{1}}^{-1}(\P_{1})) \subset \P_{2}$ and $\pi_{\B_{2}}\circ \Phi(\pi_{\B_{1}}^{-1}(\B_{1})) \subset \B_{2}$ (one could say simply that $\Phi$ preserves the product structure, or preserves fibers). The definition means that $\Phi$ induces maps $\P_{1} \to \P_{2}$ and $\B_{1} \to \B_{2}$ consistent with the given incidence relations.  
An isomorphism of incidence structures is a bijective morphism of incidence structures such that the inverse bijection is also a morphism of incidence structures. 
An automorphism of an incidence structure is an isomorphism of the incidence structure onto itself. These form a group $\Aut(\P, \B, \F)$. The maps $\P \to \P$ and $\B \to \B$ induced by $\Phi \in \Aut(\P, \B, \F)$ are bijections, so permutations; an automorphism $\Phi$ can therefore be regarded as a permutation of $\P$ that preserves incidence. 

Given an incidence structure $(\P, \B, \F)$, write $\B_{p} = \pi_{\B}(\pi_{\P}^{-1}(p))$ for the set of blocks containing $p \in \P$ and $\P_{B} = \pi_{\P}(\pi_{\B}^{-1}(B))$ for the set of points contained in $B \in \B$. Write $\card S$ for the cardinality of the set $S$. A \emph{partial linear space} is an incidence structure $(\P, \B, \F)$ such that for any two distinct $p, q \in \P$, $\card (\B_{p} \cap \B_{q}) \leq 1$. More prosaically, this means that two distinct element of $\P$ are incident with at most one element of $\B$. This implies that $\card (\P_{B}\cap \P_{C}) \leq 1$ for any distinct $B, C \in \B$, for were $\card (\P_{B}\cap \P_{C}) \geq 2$ then $B, C \in \B_{p}\cap \B_{q}$. Hence the dual incidence structure $(\B, \P, \F^{\ast})$ is also a partial linear space. A subspace of a partial linear space is necessarily again a partial linear space. A partial linear space is a \emph{linear space} if for any two distinct $p, q \in \P$, $\card (\B_{p} \cap \B_{q}) = 1$ (equivalently  $\card (\P_{B}\cap \P_{C}) = 1$ for any distinct $B, C \in \B$). More prosaically, this means that every pair of distinct elements of $\P$ is incident to a unique element of $\B$.

A \emph{partial Steiner triple system} on $\P$ is a partial linear space $(\P, \B, \F)$ in which each block has exactly three points, and a \emph{Steiner triple system} on $\P$ is a linear space $(\P, \B, \F)$ in which each block has exactly three points. 

The \emph{$m$-dimensional affine space of order $p$}, $AG(m, p)$, is the partial linear space $(\P, \B, \F)$ where $\P$ is a $m$-dimensional vector space over the finite field with $p$ elements, $\B$ is the set of affine lines in $\P$, and $\F$ is the incidence relation corresponding to containment of points in lines. When $m=2$ it is called the \emph{affine plane of order $p$}. When $p = 3$, the $AG(m, 3)$ is a Steiner triple system of order $3^{m}$. The affine plane of order $3$, $AG(2, 3)$, admits the following description. Fill a $3 \times 3$ grid with the numbers from $\bar{9}$ and tile the plane with this grid. In each row, column, and diagonal there are exactly three distinct integers, and these are the $12$ blocks 
\begin{align}\label{9124}
\B = \{123, 456, 789, 147, 258, 369, 159, 267, 348, 168, 249, 357\}
\end{align}
of a Steiner triple system $\B$ having format $(9, 12, 4)$. The \emph{projective plane of order p} is the partial linear space $(\P, \B, \F)$ where $\P$ is the set of nonzero points in a $3$-dimensional vector space over the finite field with $p$ elements, $\B$ is the set of lines in $\P$, and $\F$ is the incidence relation corresponding to containment of points in lines. The dual of the affine plane of order $2$ is isomorphic to the partial linear space obtained from the projective plane of order $2$ by deleting a point and all the lines through this point. 

Up to isomorphism the set of points $\P$ of a partial linear space is determined by its cardinality, and when speaking of triple systems on a set of $n$ points it is convenient to take as $\P$ the set $\setn = \{1, \dots, n\}$, regarded as a finite set with $n$ distinguishable elements (the ordering implicit in their labels is forgotten). The set of all (partial) Steiner triple systems on $\setn$ is denoted $STS(n)$ ($PSTS(n)$). Thus $\B \in PSTS(n)$ means a partial linear space $(\setn, \B, \F)$ such that each element of $\B$ has cardinality $3$. In this case it is convenient to identify $\B$ with a subset of the power set $2^{\setn}$ and to omit $\F$ from the notation. 

The definition of a partial Steiner triple system does not require that the number of blocks of $\B$ containing a given element of $\setn$ be independent of the chosen element. A $\B \in PSTS(n)$ for which there is an integer $r \geq 1$ such that each element of $\setn$ is contained in exactly $r$ blocks of $\B$ is said to be \emph{regular} with \emph{replication number} $r$. The subset of $PSTS(n)$ comprising regular systems with replication number $r$ is denoted $PSTS_{r}(n)$.
A \emph{Steiner triple system} on $\setn$ is a regular partial Steiner triple system $\B \in PSTS_{r}(n)$ such that each pair of distinct $i \neq j \in \setn$ is contained in exactly one block $B$ of $\B$ and each $i \in \setn$ is contained in exactly $r$ blocks of $\B$.
The number of blocks $b = \card \B$ in $\B \in STS(n)$ and its replication number $r$ satisfy $b = n(n-1)/6$ and $r = (n-1)/2$.

The set $STS(n)$ is nonempty if and only if $n$ equals $1$ or $3$ modulo $6$, while $PSTS(n)$ is nonempty for all $n \geq 0$. 
In particular, $PSTS(n)$ contains the empty partial Steiner triple system $\emptyset \subset \setn$ having no blocks. It is convenient to regard $\emptyset$ as having replication number $r = 0$.

For $\B \in PSTS(n)$ write $i \sim j$ if a pair $i \neq j \in \setn$ is contained in a block of $\B$, and $i \nsim j$ otherwise.
If $i \sim j$ then the block containing $i$ and $j$ is unique, so there is a unique element $i \join j \in \setn$ such that $\{i, j, i\join j\} \in \B$. For $\B \in STS(n)$ define $i \join i = i$. In this case the commutative operation $\join$ satisfies $i\join (i \join j) = j$, so $(\setn, \join)$ is an idempotent totally symmetric (or commutative) quasigroup, and the structure of the Steiner triple system $\B$ is completely determined by the associated quasigroup, which called here the \emph{underlying Steiner quasigroup}.

Throughout the paper, $\fie$ is a field of characteristic zero and $\fiet$ is the multiplicative group of nonzero elements of $\fie$.

\begin{lemma}\label{hatlemma}
For $\B \in PSTS(n)$ and $\ga, \al, \be \in \fie$ not all equal to zero, let $(\ualg_{\ga, \al, \be}(\fie, \B), \hmlt)$ be the commutative algebra generated by $\{\hat{e}_{i}: i \in \setn\}$ subject to the relations
\begin{align}\label{ualgrelations}
&\hat{e}_{i} \hmlt \hat{e}_{i} = \ga\hat{e}_{i}, &
&\hat{e}_{i} \hmlt \hat{e}_{j} = \begin{cases}
\al(\hat{e}_{i} + \hat{e}_{j}) + \be \hat{e}_{i\join j}&  i \sim j, \,\, i \neq j \in \setn,\\
0 & i \nsim j , \,\, i \neq j \in \setn.
\end{cases}
\end{align}
\begin{enumerate}
\item\label{idealclaim} The span of $\hat{e} = \sum_{i \in \setn}\hat{e}_{i}$ is an ideal in $(\ualg_{\al, \be}(\fie, \B), \hmlt)$ if and only if there holds one of the following two possibilities.
\begin{enumerate}
\item\label{idealclaim1} $\B \in STS(n)$ and $\be = \ga + (n-2)\al$. In this case $x\mlt \hat{e} = (\al + \be)(\sum_{i \in \setn}x_{i})\hat{e}$ for all $x = \sum_{i \in \setn}x_{i}\hat{e}_{i} \in \ualg_{\ga, \al, \be}(\fie, \B)$.
\item \label{idealclaim2} $\B$ is regular with replication number $r$, $\be = -\al$, and $\ga = -2r\al$. In this case $x\mlt \hat{e} = 0$ for all $x \in \ualg_{-2r\al, \al, -\al}(\fie, \B)$. 
\end{enumerate}
\item\label{hatexactclaim} $(\ualg_{\ga, \al, \be}(\fie, \B), \mlt)$ is exact if and only if either $\B$ is regular with replication number $r$ satisfying $\ga + 2r\al = 0$ or $\ga = \al = 0$.
\end{enumerate}
\end{lemma}
\begin{proof}
For fixed $i$,
\begin{align}\label{hehehe}
\sum_{j \neq i, j \sim i}\hat{e}_{j} = \sum_{j \neq i, j\sim i}\hat{e}_{i\join j} = \hat{e} - \hat{e}_{i} - \sum_{j \neq i, j \nsim i}\hat{e}_{j}.
\end{align}
Let $r_{i}$ be the number of blocks of $\B$ containing $i$. By \eqref{hehehe},
\begin{align}\label{heideal}
\begin{split}
\hat{e}_{i}\hmlt \hat{e} & = \ga \hat{e}_{i} + \sum_{j \neq i, j \sim i}\left(\al(\hat{e}_{i} + \hat{e}_{j}) + \be \hat{e}_{i\join j}\right) = (\ga + 2r_{i}\al)\hat{e}_{i} + \al\sum_{j \neq i, j \sim i}\hat{e}_{j}  + \be \sum_{j \neq i, j\sim i}\hat{e}_{i\join j}\\
& = (\ga + 2r_{i}\al)\hat{e}_{i} + (\al + \be)(\hat{e} - \hat{e}_{i} - \sum_{j \neq i, j \nsim i}\hat{e}_{j})\\
& = (\al + \be)\hat{e} + (\ga + (2r_{i}-1)\al - \be)\hat{e}_{i} - ( \al + \be) \sum_{j \neq i, j \nsim i}\hat{e}_{j}.
\end{split}
\end{align}
By \eqref{heideal}, $\fie \hat{e}$ is an ideal if and only if there hold 
\begin{align}
\label{hathe1}&\ga + (2r_{i}-1)\al - \be = 0,& \\
\label{hathe2}& ( \al + \be) \sum_{j \neq i, j \nsim i}\hat{e}_{j} = 0,
\end{align}
for all $1 \leq i \leq n$. Suppose $\al \neq 0$. Then \eqref{hathe1} holds for all $i$ if and only if $r_{i} = r = (\be - \ga)/(2r - 1)$, in which case $\B$ is regular with replication number $r$. If $\be \neq -\al$ then \eqref{hathe2} can only be true if $j \sim i$ for all $j \neq i$, that is $r = (n-1)/2$ and $\B$ is a Steiner triple system. In this case there holds \eqref{idealclaim1} with $\al \neq 0$. If $\be = -\al$, then there holds \eqref{idealclaim2}. Now suppose $\al = 0$. Then \eqref{hathe1} holds for all $i$ if and only if $\be = \ga$. Since by assumption not all of $\al$, $\be$, and $\ga$ are zero, this means $\be \neq 0$, so \eqref{hathe2} holds if and only if  $j \sim i$ for all $j \neq i$, that is $r = (n-1)/2$ and $\B$ is a Steiner triple system. In this case there holds \eqref{idealclaim1} with $\al = 0$.
Because $\tr L(\hat{e}_{i}) = \ga + 2r_{i}\al$, $(\ualg_{\ga, \al, \be}(\fie, \B), \mlt)$ is exact if and only if $\ga + 2r_{i}\al = 0$ for all $i$. If $\al \neq 0$ this implies $r_{i} = r$ does not depend on $i$, and $\B$ is regular with replication number $r$. This shows \eqref{hatexactclaim}.
\end{proof}

\begin{remark}
For $r \in \fiet$ the map $\Psi_{r}:(\ualg_{\ga, \al, \be}(\fie, \B), \mlt) \to (\ualg_{r\ga, r\al, r\be}(\fie, \B), \mlt)$ defined by $\Psi(\hat{e}_{i}) = r\hat{e}_{i}$ is a linear isomorphism. Consequently, any one of the parameters $\ga$, $\al$, or $\be$ can be normalized to take values in $\{0, 1\}$. It is most convenient to normalize $\ga$ in this way, in which case each $\hat{e}_{i}$ is an idempotent if $\ga = 1$ or a square-zero element if $\ga = 0$.
\end{remark}

\begin{definition}\label{talgdefined}
For $\be \in \fie$ and $\B \in STS(n)$ or for $\be = 0$ and $\B = \emptyset \in PSTS_{0}(n)$, define $\ualg_{\be}(\fie, \B) = \ualg_{1, \al, \be}(\fie, \B)$ where $\al = \tfrac{\be - 1}{n-2}$ and define $\talg_{\be}(\fie, \B)$ to be the quotient of $\ualg_{\be}(\fie, \B)$ by $\fie\hat{e}$, equipped with the multiplication $\mlt$ induced from that on $\ualg_{\be}(\fie, \B)$. 
\end{definition}

Lemma \ref{eigenlemma} motivates using the parameters $\be_{\pm}$ defined by
\begin{align}\label{bepm}
&\be_{+} = -\tfrac{1}{n-2}  + \tfrac{n-1}{n-2}\be = \al + \be,& &\be_{-} = -\tfrac{1}{n-2} - \tfrac{n-3}{n-2}\be  = \al - \be = -\tfrac{2}{n-1} - \tfrac{n-3}{n-1}\be_{+}.
\end{align}
Any one of the parameters $\al$, $\be$, and $\be_{\pm}$ determines the other three, and in what follows the parameters $\al$, $\be$, and $\be_{\pm}$ are used as is convenient.
The images $e_{i} \in \talg_{\be}(\fie, \B)$ of the $\hat{e}_{i}$ satisfy the relations
\begin{align}\label{talgrelations}
\begin{aligned}
& \sum_{i \in \setn} e_{i} = 0, \qquad e_{i} \mlt e_{i} = e_{i},& \\
&\begin{aligned}e_{i} \mlt e_{j} &=  
\al(e_{i} + e_{j}) + \be e_{i\join j} = -\tfrac{1}{n-2}(e_{i} + e_{j}) + \tfrac{\be}{n-2}(e_{i} + e_{j} +  (n-2)e_{i\join j})\\
&= \tfrac{\be_{+}}{2}(e_{i} + e_{j} + e_{i \join j}) + \tfrac{\be_{-}}{2}(e_{i} + e_{j} - e_{i\join j}),& &i \neq j \in \setn,\end{aligned}
\end{aligned}
\end{align}
and any $n-1$ of the $e_{i}$ constitute a basis of $\talg_{\be}(\fie, \B)$. 

\begin{remark}\label{griessharadaremark}
For any $n \geq 2$, taking $\be = 0$ in \eqref{talgrelations} yields the $(n-1)$-dimensional algebra $(\talg_{0}(\emptyset), \mlt)$ defined by the relations
\begin{align}\label{ealgrelations}
& \sum_{i \in \setn} e_{i} = 0,& &e_{i} \mlt e_{i} = e_{i},&&e_{i} \mlt e_{j} = -\tfrac{1}{n-2}(e_{i} + e_{j}).
\end{align}
Here this algebra is denoted $\ealg^{n-1}(\fie)$ and is called the \emph{simplicial} algebra. It was first studied by K. Harada \cite{Harada, Harada-doublytransitive, Harada-automorphisms}, and R.~L. Griess in unpublished work (see the second paragraph of \cite{Smith} and \cite[Appendix]{Dong-Griess}). The proofs of Theorems \ref{stsalgebratheorem} and \ref{stsalgebrasimpletheorem} apply (as written) to $\ealg^{n-1}(\fie)$ for all $n$, recovering its simplicity and other basic features first shown by Harada and Griess.

It is apparent from \eqref{talgrelations} that $(\talg_{\be}(\fie, \B), \mlt)$ can be viewed as a deformation of $\ealg^{n-1}(\fie) = (\talg_{0}(\fie, \B), \mlt) $.

The unique Steiner triple system on $\setn$ for $n < 7$ is that on $\bar{3}$ containing the single block $\{1,2,3\}$. In this case $e_{i} + e_{j} + e_{i\join j} =0$ if $i \neq j \in \bar{3}$, so \eqref{talgrelations} becomes $e_{i}\mlt e_{j} = -e_{i} - e_{j}$, showing that $(\talg_{\be}(\{1, 2, 3\}), \mlt) = (\talg_{0}(\{1, 2, 3\}), \mlt) = \ealg^{2}(\fie)$ for all $\be \in \fie$. 
\end{remark}

\begin{remark}\label{mendelsohnremark}
In \cite{Mendelsohn} E. Mendelsohn associated with $\B \in STS(n)$ a commutative nonassociative algebra. Although its definition resembles that of $(\ualg_{0,1}, \mlt)$, they are not the same.
Taking $\ga = 0$, $\al = 0$, and $\be = 1$ in \eqref{ualgrelations} yields an $n$-dimensional algebra defined by the relations 
\begin{align}\label{mendelrelations}
&\hat{e}_{i} \mlt \hat{e}_{i} =0,& &\hat{e}_{i} \mlt \hat{e}_{j} =  \hat{e}_{i\join j},& &i \neq j \in \setn.
\end{align}
Note that the $\hat{e}_{i}$ are square-zero elements rather than idempotents. Mendelsohn's algebra is the $(n+1)$-dimensional algebra defined by adding a unit to the algebra \eqref{mendelrelations}.
\end{remark}

\begin{remark}
Taking $\be = -\al$ and $\ep = 1$ in \eqref{ualgrelations} yields the $n$-dimensional algebra defined by the relations
\begin{align}\label{matsuorelations}
&\hat{e}_{i} \mlt \hat{e}_{i} = \hat{e}_{i},& &\hat{e}_{i} \mlt \hat{e}_{j} = \al(e_{i} + e_{j} -  \hat{e}_{i\join j}),& &i \neq j \in \setn.
\end{align}
This is the \emph{Matsuo algebra} $\matsuo_{\al}(\fie, \B)$ of $\B \in STS(n)$ defined in \cite[section $3.2$]{Matsuo-3transpositionarxiv} (different authors normalize differently the constant in the definition of $\matsuo_{\al}(\fie, \B)$). See \cite[Section $6$]{Hall-Rehren-Shpectorov-primitive}, \cite{DeMedts-Rehren}, and \cite[Section $5$]{Krasnov-Tkachev} for details and more information.
The proof of Lemma \ref{hatlemma} shows that $\hat{e}\hmlt \hat{e}_{i} = (1 - (n-1)\be)\hat{e}_{i}$ for $i \in \setn$. In particular $\matsuo_{\al}(\fie, \B)$ is unital except when $\al = -\be = -1/(n-1)$ (this is part of \cite[Proposition $3.2.1$]{Matsuo-3transpositionarxiv}). On the other hand, it is uniquely in this exceptional case $\al = -\be = -1/(n-1)$ that $\fie\hat{e}$ is an ideal in the Matsuo algebra.
 
These observations give some motivation for regarding $(\talg_{1/(n-1)}(\fie, \B), \mlt)$ as meriting special consideration. Another reason for such special consideration is given by Lemma \ref{eigenlemma}.
\end{remark}

The automorphism group $\Aut(\B) = \Aut(\setn, \B)$ of $\B \in PSTS(n)$ is its automorphism group viewed as a partial linear system, so comprises those permutations of $\setn$ that map blocks of $\B$ to blocks of $\B$. Given an element $\si \in \Aut(\B)$ (viewed as a permutation of $\setn$ preserving the blocks of $\B$) define $\hat{\phi}_{\si} \in \eno(\ualg_{\be})$ by $\hat{\phi}(\hat{e}_{i}) = \hat{e}_{\si(i)}$ and extending linearly. Since $\si(i \join j) = \si(i)\join \si(j)$, $\hat{\phi}_{\si}$ is an algebra automorphism of $(\ualg_{\be}, \mlt)$. Since $\hat{\phi}_{\si}(\hat{e}) = \hat{e}$, $\hat{\phi}_{\si}$ descends to an algebra automorphism $\phi_{\si}$ of $(\talg_{\be}(\fie, \B), \mlt)$ satisfying $\phi_{\si}(e_{i}) = e_{\si(i)}$ for $1 \leq i \leq n$. This shows that the automorphism group $\Aut(\B)$ of the Steiner triple system $\B$ acts on $(\talg_{\be}(\fie, \B), \mlt)$ by algebra automorphisms, meaning that the map $\si \to \phi_{\si}$ is an embedding $\Aut(\B) \to \Aut(\talg_{\be}(\fie, \B), \mlt)$. For what $\be$ and $\B$ this map is an isomorphism has not been resolved. Lemma \ref{prepermutelemma} characterizes its image.
\begin{lemma}\label{prepermutelemma}
Let $\be \in \fie\setminus\{0\}$ and consider $\B \in STS(n)$. If $\phi \in \Aut(\talg_{\be}(\fie, \B), \mlt)$ permutes $\{e_{i}: i \in \setn\}$, then the permutation $\si \in S_{n}$ defined by $\phi(e_{i}) = e_{\si(i)}$ is an automorphism of $\B$.
\end{lemma}
\begin{proof}
For all $i \neq j \in \setn$, $0 = \phi(e_{i}\mlt e_{j}) - \phi(e_{i})\mlt \phi(e_{j}) = \be(e_{\si(i\join j)} - e_{\si(i)\join \si(j)})$. Because $\be \neq 0$, $\si(i \join j) = \si(i)\join \si(j)$, which means that $\si$ is an automorphism of $\B$.
\end{proof}

\begin{remark}
An embedding of $\C \in STS(m)$ in $\B \in STS(n)$ is an injection $\imt: \setm \to \setn$ such that $\imt(\C) \subset \B$. It determines an algebra embedding $\talg_{\be}(\fie, \C) \to \talg_{\be}(\fie, \B)$ defined as the extension by linearity of the map defined by $e_{i} \to e_{\imt(i)}$ for $i \in \setm$. 
\end{remark}

\section{Killing metrizability and simplicity}

Lemma \ref{eigenlemma} is similar to \cite[Theorem $6.2$]{Hall-Rehren-Shpectorov-primitive}. Recall the definition of $\be_{\pm}$ in \eqref{bepm}.
\begin{lemma}\label{eigenlemma}
Let $\be \in \fie$ and $\B \in STS(n)$. 
For $e_{i} \in (\talg_{\be}(\fie, \B), \mlt, \kform)$ define the $(n-1)/2$-dimensional subspace 
\begin{align}\label{alminusbeeigen}
\balg^{e_{i}}_{-} = \fie\{e_{j} - e_{i\join j}: j \neq i\},
\end{align}
and the $(n-3)/2$-dimensional subspace 
\begin{align}\label{alplusbeeigen}
\balg^{e_{i}}_{+} = \{\sum_{i \in B \in \B}\la_{B}\ga_{\B}:\sum_{i \in B \in \B}\la_{B} = 0\}= \fie\{2e_{i} + (n-1)(e_{j} + e_{i\join j}): j \neq i\}.
\end{align}
If $\be \in \fie \setminus \{-\tfrac{n-1}{n-3}, 0, 1\}$ then $1$, $\be_{+}$, and $\be_{-}$ are pairwise distinct and $\balg^{e_{i}}_{+}$ and $\balg^{e_{i}}_{-}$ are the $\be_{+}$ and $\be_{-}$ eigenspaces of $L(e_{i})$. In the excluded cases:
\begin{enumerate}
\item If $\be = -\tfrac{n-1}{n-3}$, then $\be_{-}= 1$, $\balg^{e_{i}}_{-}$ has codimension $1$ in the $1$-eigenspace of $L(e_{i})$, and $\balg^{e_{i}}_{+}$ is the $\be_{+} = -(n+1)/(n+3)$ eigenspace of $L(e_{i})$.
\item If $\be = 0$ then $\be_{+} = \be_{-} = -1/(n-2)$ and $\balg^{e_{i}}_{-}\oplus \balg^{e_{i}}_{+}$ is the $-1/(n-2)$ eigenspace of $L(e_{i})$.
\item If $\be = 1$, then $\be_{+} = 1$, $\balg^{e_{i}}_{+}$ has codimension $1$ in the $1$-eigenspace of $L(e_{i})$ and $\balg^{e_{i}}_{-}$ is the $\be_{-} = -1$ eigenspace of $L(e_{i})$.
\end{enumerate}
\end{lemma}

\begin{proof}
The cardinality $(n-1)/2 = r$ set $\{e_{j} - e_{i \join j}: i \in B \,\,\text{and}\,\, j < i\join j\}$ is linearly independent and spans $\balg^{e_{i}}_{-}$. Computations using \eqref{talgrelations} and \eqref{bepm} show $e_{i}\mlt(e_{j} - e_{i\join j}) =\be_{-}(e_{j} - e_{i\join j})$. This shows $\balg^{e_{i}}_{-}$ is contained in $\ker(L(e_{i}) - \be_{-}I)$. 

For $B = \{i, j, i\join j\} \in \B$, the element $\ga_{B} = e_{i} + e_{j} + e_{i\join j} \in \talg_{\be}(\fie, \B)$ satisfies 
\begin{align}\label{gabgab}
\ga_{B}\mlt \ga_{B} 
= (1 +3\be_{+} + \be_{-})\ga_{B} = \tfrac{2n\be + n-6}{2n}\ga_{B}.
\end{align}There holds $\sum_{i \in B \in \B}\ga_{B} = re_{i} + \sum_{j \neq i}e_{j} = (r-1)e_{i} = \tfrac{n-3}{2}e_{i}$. Consequently, for $B = \{i, j, i\join j\} \in \B$, $2\sum_{i \in B^{\prime} \neq B}(\ga_{B} - \ga_{B^{\prime}}) = 2r\ga_{B} - 2(r-1)e_{i} = 2e_{i} + (n-1)(e_{j} + e_{i \join j})$. This shows that the left-hand side of \eqref{alplusbeeigen} is contained in the right-hand side; the reverse inclusion is clear. For $k \neq i$, the cardinality $(n-3)/2 = r - 1$ set $\{2e_{i} + (n-1)(e_{j} + e_{i\join j}): i \in B, k \notin B\}$ is linearly independent and spans $\balg^{e_{i}}_{+}$. For $i =  B \cap B^{\prime}$, it follows from $e_{i}\mlt \ga_{B} = (1 + \tfrac{1}{2}\be_{+} + \tfrac{3}{2}\be_{-})e_{i} + \be_{+}\ga_{B}$ that $e_{i}\mlt(\ga_{B} - \ga_{B^{\prime}}) = \be_{+}\ga_{B}$. This shows $\balg^{e_{i}}_{+}$ is contained in the $\be_{+}$ eigenspace of $L(e_{i})$.

The preceding shows that the eigenvalues of $L(e_{i})$ are $1$, $\be_{-}$, and $\be_{+}$, with respective multiplicities at least $1$, $(n-1)/2$, and $(n-3)/2$, with equality of the multiplicities except when the values $\{1, \be_{+}, \be_{-}\}$ are not pairwise distinct, which occurs in the indicated exceptional cases.
\end{proof}

\begin{remark}
Lemma \ref{eigenlemma} shows that in $(\talg_{\be}(\fie, \B), \mlt)$ the multiplication operator $L(e_{i})$ of the idempotent $e_{i}$ is semisimple with minimal polynomial dividing $(x -1)(x - \be_{+})(x - \be_{-})$. As is indicated in Table \ref{stseigentable2}, qualitative transitions occur when one of $\be_{\pm}$ is in the set $\{0, 1/2, 1\}$ or they are equal. Theorem \ref{stsalgebrasimpletheorem} shows that $(\talg_{\be}(\fie, \B), \mlt)$ is simple except when $\be_{+}$ or $\be_{-}$ equals $1$. 
\begin{table}

\begin{tabular}{|c||c|c|c|c|c|c|c|}
\hline
$\be$ &  $-\frac{n-1}{n-3}$ &$-\frac{n}{2(n-3)}$ & $-\frac{1}{n-3}$ & $0$ & $\frac{1}{n-1}$ & $\tfrac{n}{2(n-1)}$ &$1$\\
\hline
$\be_{+} = - \tfrac{1}{n-2} + \tfrac{n-1}{n-2}\be$ &  $-\frac{n+1}{n-3}$ & $-\frac{n+3}{2(n-3)}$ & $-\frac{2}{n-3}$ & $-\tfrac{1}{n-2}$ & $0$ & $\tfrac{1}{2}$ & $1$\\
\hline
$\be_{-} = - \tfrac{1}{n-2} - \tfrac{n-3}{n-2}\be $ & $1$ &  $\tfrac{1}{2}$ & $0$ &$-\tfrac{1}{n-2}$ & $-\frac{2}{n-1}$ & $-\tfrac{n+1}{2(n-1)}$& $-1$\\
\hline
\end{tabular}
\caption{Transitional values of eigenvalues of $L(e_{i})$.}\label{stseigentable2}
\end{table}
\end{remark}

\begin{corollary}\label{exactcorollary}
For $\be \in \fie$ and $\B \in STS(n)$, the algebra $(\talg_{\be}(\fie, \B), \mlt)$ is exact.
\end{corollary}

\begin{proof}
By \eqref{bepm}, $(n-3)\be_{+} + (n-1)\be_{-} = -2$, so, by Lemma \ref{eigenlemma}, $\tr L(e_{i}) = 1 + \tfrac{n-1}{2}\be_{-} + \tfrac{n-3}{2}\be_{+} = 0$. Since the $e_{i}$ span $\talg_{\be}$, this suffices to show that $(\talg_{\be}(\fie, \B), \mlt)$ is exact.
\end{proof}

\begin{theorem}\label{stsalgebratheorem}
For $\be \in \fie$ and $\B \in STS(n)$, the Killing form $\kform(x, y) = \tr L(x)L(y)$ of the algebra $(\talg_{\be}(\fie, \B), \mlt)$ is invariant. If $\om = (n-3)\be^{2} + 1 \neq 0$, then $\kform$ is nondegenerate. If $\fie$ is a Euclidean field, then $\kform$ is positive definite.
\end{theorem}
\begin{proof}
By Lemma \ref{eigenlemma}, 
\begin{align}
\label{ststauii}
\begin{split}
\kform(e_{i}, e_{i}) &= \tr L(e_{i})^{2} = 1 + \tfrac{n-1}{2}\be_{-}^{2} + \tfrac{n-3}{2}\be_{+}^{2}= \tfrac{n-1}{n-2}\left((n-3)\be^{2} + 1\right) =  \tfrac{n-1}{n-2}\om.
\end{split}
\end{align}
The trace $\tr L(e_{i})L(e_{j})$ can be calculated as the trace of the matrix representing $L(e_{i})L(e_{j})$ with respect to a basis of $\talg_{\be}(\fie, \B)$. For this purpose it is convenient to use the basis $\{e_{k}: 1 \leq k \leq n-1\}$. Write $\al = \tfrac{\be - 1}{n-2}$, so that $\be = 1 + (n-2)\al$. For $i \neq j \in \setn$ and $k \notin \{i, j, i\join j\}$, straightforward calculations using \eqref{talgrelations} show
\begin{align}
\label{stsij}
\begin{split}
L(e_{i})L(e_{j})e_{i} & = \al(\al + \be + 1)e_{i} + (\al^{2} + \be^{2})e_{j} + 2\al \be e_{i\join j},\\
L(e_{i})L(e_{j})e_{j} & = \al e_{i} + \al e_{j} + \be e_{i\join j},\\
L(e_{i})L(e_{j})e_{i\join j} & = (2\al^{2} + \be)e_{i} + \al(\al + \be)e_{j} + \al(\al + \be)e_{i\join j},
\end{split}\\
\label{stsijk}
\begin{split}
L(e_{i})L(e_{j})e_{k} & = \al(2\al + \be)e_{i} + \al^{2}e_{j} + \al^{2}e_{k} + \al\be(e_{i\join j} + e_{i\join k} + e_{j\join k}) + \be^{2}e_{i\join(j\join k)},
\end{split}
\end{align}
In particular, $L(e_{i})L(e_{j})$ preserves $\fie\{e_{i}, e_{j}, e_{i\join j}\}$. In using these computations to compute the trace of $L(e_{i})L(e_{j})$, care needs to be taken because $i \join j$, $i\join k$, $j \join k$, $i\join(j\join k)$ can equal $n$, in which case the corresponding vectors need to be rewritten in terms of the basis $\{e_{k}: 1 \leq k \leq n-1\}$. Note that $i \join j = n$ if and only if $j = i \join n$, while $i\join k = n$, $j \join k = n$, and $i\join(j\join k) = n$ hold if and only if $k = i\join n$, $k = j\join n$, and $k = j\join(i\join n)$, respectively.
In these cases, for $1 \leq i\neq j \leq n-1$, specializing \eqref{stsij} and \eqref{stsijk} yields
\begin{align}
\label{stsijn}
\begin{split}
L(e_{i})L(e_{i\join n})e_{i} & = \al(\al + \be + 1)e_{i} + (\al^{2} + \be^{2})e_{i\join n} + 2\al \be e_{n} \\&= \al(\al -\be + 1)e_{i} +(\al-\be)^{2}e_{i\join n} -2\al\be\sum_{l \notin\{i, n, i\join n\}}e_{l} ,\\
L(e_{i})L(e_{i\join n})e_{i\join n} & = \al e_{i} + \al e_{i\join n} + \be e_{n} = (\al - \be)(e_{i} + e_{i\join n}) -\be\sum_{l \notin\{i, n, i\join n\}}e_{l},
\end{split}
\end{align}
and, when $j \neq i \join n$,
\begin{align}
\label{stsijkn}
\begin{split}
L(e_{i})&L(e_{j})e_{i\join n} \\& = \al(2\al + \be)e_{i} + \al^{2}e_{j} + \al^{2}e_{i\join n} + \al\be(e_{i\join j} + e_{n} + e_{j\join (i \join n)}) + \be^{2}e_{i\join(j\join (i\join n))}\\
& = \al(\al-\be)e_{i\join n} + \text{terms involving neither $e_{i\join n}$ nor $e_{n}$},\\
L(e_{i})&L(e_{j})e_{j\join n} \\& = \al(2\al + \be)e_{i} + \al^{2}e_{j} + \al^{2}e_{j\join n} + \al\be(e_{i\join j} + e_{i\join (j\join n)} + e_{n}) + \be^{2}e_{i\join n}\\
& = \al(\al-\be)e_{j\join n} + \text{terms involving neither $e_{j\join n}$ nor $e_{n}$},\\
L(e_{i})&L(e_{j})e_{j\join (i\join n)} \\& = \al(2\al + \be)e_{i} + \al^{2}e_{j} + \al^{2}e_{j\join (i\join n)} + \al\be(e_{i\join j} + e_{i\join (j \join(i\join n))} + e_{i\join n}) + \be^{2}e_{n}\\
& = (\al^{2}-\be^{2})e_{j\join(i\join n)} + \text{terms involving neither $e_{j\join (i\join n)}$ nor $e_{n}$},
\end{split}
\end{align}
and, when $k \notin\{i, i\join n, n\}$,
\begin{align}\label{stsiink}
\begin{split}
L(e_{i})&L(e_{i\join n})e_{k}  \\
&= \al(2\al + \be)e_{i} + \al^{2}e_{i\join n} + \al^{2}e_{k} + \al\be(e_{n} + e_{i\join k} + e_{(i\join n)\join k}) + \be^{2}e_{i\join((i\join n)\join k)}\\
& = \al(\al - \be)e_{k} + \text{terms involving neither $e_{k}$ nor $e_{n}$}.
\end{split}
\end{align}
Suppose $1 \leq i \neq j \leq n-1$ and $i \join j \neq n$, so that $j \neq i \join n$. From \eqref{stsij}, \eqref{stsijk}, and \eqref{stsijkn}, it follows that the diagonal elements of the matrix of $L(e_{i})L(e_{j})$ in the basis $\{e_{l}: 1 \leq l \leq n-1\}$ are
$\al(\al + \be +1 )$ (for index $i$), $\al$ (for index $j$), $\al(\al + \be)$ (for index $i\join j$), $\al^{2}$ (for the $n-7$ choices of $k \notin \{i, j, i\join j, i\join n, j \join n, j \join (i\join n)\}$), $\al(\al - \be)$ (for indices $i\join n$ and $j \join n$), and $\al^{2} - \be^{2}$ (for index $j\join (i \join n)$). Summing yields
\begin{align}\label{ststauij}
\begin{split}
\kform(e_{i}, e_{j}) & = \al(\al + \be + 1) + \al + \al(\al + \be) + (n-7)\al^{2} + 2\al(\al - \be) + \al^{2}  - \be^{2} \\
&= (n-2)\al^{2} + 2\al - \be^{2} = -\tfrac{1}{n-2}\om.
\end{split}
\end{align}
Suppose $1 \leq i \leq n-1$ and $j = i \join n$, so that $i \join j = n$. From \eqref{stsijn}, \eqref{stsijkn}, and \eqref{stsiink}, it follows that the diagonal elements of the matrix of $L(e_{i})L(e_{i\join n})$ in the basis $\{e_{l}: 1 \leq l \leq n-1\}$ are
$\al(\al - \be +1 )$ (for index $i$), $\al - \be$ (for index $i\join n$), and $\al(\al - \be)$ (for the $n-3$ choices of $k \notin \{i, i\join n, n\}$). Summing yields
\begin{align}\label{ststauiin}
\begin{split}
\kform(e_{i}, e_{i\join n}) & = \al(\al - \be + 1) + \al - \be + (n-3)\al(\al  -\be) \\&= ((n-2)\al + 1)(\al - \be) + \al = -\tfrac{1}{n-2}\left((n-3)\be^{2} + 1\right).
\end{split}
\end{align}
Together \eqref{ststauii}, \eqref{ststauij}, and \eqref{ststauiin} show that the Gram matrix of the basis $\{e_{l}: 1 \leq l \leq n-1\}$ with respect to $\kform$ is $\om(n\imat_{n-1} - \ones_{n-1})$, where $\imat_{k}$ is the $k \times k$ identity matrix and $\ones_{k}$ is the $k \times k$ matrix with all entries equal to $1$. The eigenvalues of $n\imat_{n-1} - \ones_{n-1}$ are $1$ with multiplicity $1$ and $n$ with multiplicity $n-2$, so $n\imat_{n-1} - \ones_{n-1}$ is nondegenerate. This shows that $\kform$ is nondegenerate provided $\om \neq 0$. If $\fie$ is Euclidean, $\om = (n-3)\be^{2} + 1 \geq 1$, so $\kform$ is positive definite.

Using \eqref{ststauii}, \eqref{ststauij}, and \eqref{ststauiin} it can be checked that, for $1 \leq i \leq n-1$,
\begin{align}
&\kform(e_{n}, e_{n}) = \tfrac{n-1}{n-2}\om,& & \kform(e_{i}, e_{n}) = -\tfrac{1}{n-2}\om,
\end{align}
so that the Gram matrix of $\{e_{i}:1 \leq i \leq n\}$ with respect to $\kform$ is $\om(n\imat_{n} - \ones_{n})$. The matrix $n\imat_{n} - \ones_{n}$ has eigenvalues $0$ with multiplicity $1$ and $n$ with multiplicity $n-1$. The one-dimensional radical corresponds with the relation $\sum_{i \in \setn}e_{i} = 0$.

To check the invariance of $\kform$ it suffices to check it on the spanning set $\{e_{i}: i \in \setn\}$. The preceding can be used to check that
\begin{align}
\begin{split}
\kform(e_{i}, e_{i} \mlt e_{j}) &= -\tfrac{1}{n-2}\om = \kform(e_{i}, e_{j}) = \kform(e_{i}\mlt e_{i}, e_{j}),\\
\kform(e_{i}, e_{j}\mlt e_{i \join j}) & = \tfrac{(n(n-3)\be + 2)}{(n-2)^{2}}\om = \kform(e_{i}\mlt e_{j}, e_{i\join j}),\\
\kform(e_{i}, e_{j}\mlt e_{k}) & = \tfrac{(1 - (n-1)\be)}{(n-2)^{2}}\om = \kform(e_{i}\mlt e_{j}, e_{k}),
\end{split}
\end{align}
for all $i, j, k \in \setn$ such that $i \neq j$ and $k \notin \{i, j, i\join j\}$.  This suffices to show the invariance of $\kform$.
\end{proof}

\begin{corollary}\label{framecorollary}
For $\be \in \fie$ and $\B \in STS(n)$, any $x \in \talg_{\be}(\fie, \B)$ satisfies $\sum_{i \in \setn}\kform(x, e_{i})e_{i} = \frac{n((n-3)\be^{2} + 1)}{n-2}x$.
\end{corollary}

\begin{proof}
This follows straightforwardly from the proof of Theorem \ref{stsalgebratheorem}.
\end{proof}

\begin{remark}
A set $\F$ of vectors in a real or complex Hilbert space $(\ste, \lb \dum, \dum \ra)$ is called a \emph{tight frame} with \emph{frame constant} $A$ if it satisfies the identity $\sum_{v \in \F}\lb x, v\ra v = Ax$ for some constant $A$ \cite{Waldron}. Corollary \ref{framecorollary} shows that for $\B \in STS(n)$ the generating idempotents $e_{i}$ constitute a tight frame with frame constant $\frac{n((n-3)\be^{2} + 1)}{n-2}$.
\end{remark}

\begin{lemma}\label{simplelemma}
Let $(\alg, \mlt)$ be a finite-dimensional commutative algebra over a field $\fie$ of characteristic $0$ admitting a nondegenerate symmetric bilinear form $\mu$ that is invariant, meaning $\mu(x\mlt y, z) = \mu(x, y \mlt z)$. Suppose $\axes \subset \alg$ spans $\alg$ and that each $a \in \axes$ is a $\mu$-anisotropic idempotent satisfying that $L(a)$ is diagonalizable over $\fie$ with $1$ as an eigenvalue of multiplicity $1$. Then every nontrivial ideal in $(\alg, \mlt)$ contains an element of $\axes$. If, moreover, $L(a)$ is invertible for every $a \in \axes$, then $(\alg, \mlt)$ is simple.
\end{lemma}

\begin{proof}
Let $\ideal \subset \talg_{\be}(\fie, \B)$ be a nontrivial ideal and suppose $0 \neq x \in \ideal$. By the nondegeneracy of $\mu$, because $\axes$ spans $\alg$, there is $a \in \axes$ such that $\mu(x, a) \neq 0$. 
Let $1, \la_{1}, \dots, \la_{k}$ be the distinct eigenvalues of $L(a)$, so that $(L(a) - \imat)Q = 0$ where $Q = \prod_{i}(L(a) - \la_{i}\imat)$. By hypothesis $a$ spans $\ker(L(a) - \imat)$, so $Qx = \rho a$ for some $\rho \in \fie$. By the invariance of $\mu$, $Q$ is $\mu$-self-adjoint, so $\rho \mu(a, a) = \mu(Qx, a) = \mu(x, Qa) = \left(\prod_{i}(1-\la_{i})^{-1}\right)\mu(x, a)\neq 0$. This shows $\rho \neq 0$, so $a = \rho^{-1}Qx \in \ideal$. If $L(a)$ is invertible, then for any $x \in \alg$ there is $y \in \alg$ such that $x = L(a)y \in \ideal$, which shows $\ideal = \alg$. 
\end{proof}

\begin{theorem}\label{stsalgebrasimpletheorem}
For $\be \in \fie\setminus\{1,  -\tfrac{n-1}{n-3}\}$ such that $(n-3)\be^{2} + 1 \neq 0$ and $\B \in STS(n)$, the algebra $(\talg_{\be}(\fie, \B), \mlt)$ is simple.
\end{theorem}

\begin{proof}
If $n = 3$, then $(\talg_{\be}(\fie, \B), \mlt)$ equals $(\talg_{0}(\fie, \B), \mlt) = \ealg^{3}(\fie)$, and so in this case it can be assumed throughout the proof that $\be = 0$. 

By Lemma \ref{eigenlemma}, for $i \in \setn$, $L(e_{i})$ is diagonalizable with $1$ an eigenvalue of multiplicity $1$ provided $\be \in \fie\setminus\{1,  -\tfrac{n-1}{n-3}\}$, and by Theorem \ref{stsalgebratheorem}, $\kform$ is nondegenerate and invariant, so the set $\axes = \{e_{i}: i \in \setn\}$ and form $\kform$ satisfy the hypotheses of Lemma \ref{simplelemma}, and so a nontrivial ideal $\ideal \subset \talg_{\be}(\fie, \B)$ contains $e_{i}$ for some $i \in \setn$.
By Lemma \eqref{eigenlemma}, $L(e_{i})$ is invertible if $\be_{+}\be_{-} \neq 0$, so in these cases $\ideal = \talg_{\be}(\fie, \B)$. The excluded cases occur when $\be_{-} = 0$, so that $\be = -\tfrac{1}{n-3}$, and $\be_{+} = 0$, so that $\be = \tfrac{1}{n-1}$. These cases require individual arguments.

Suppose $\be = \tfrac{1}{n-1}$. Because $e_{i} \in \ideal$, the $\be_{-} = -2/(n-1)$ eigenspace of $L(e_{i})$ is contained in $\ideal$, so $e_{j} - e_{i\join j} \in \ideal$ for all $j \neq i$. Hence $(e_{j}-e_{i\join j})\mlt (e_{j} - e_{i\join j}) \in \ideal$ too. There results
\begin{align}
\begin{split}
2(n+1)e_{i\join j} & = (n-1)(e_{j}-e_{i\join j})\mlt (e_{j} - e_{i\join j})  + 2e_{i} - (n+1)(e_{j} - e_{i\join j}) \in \ideal.
\end{split}
\end{align}
Hence $e_{j} \in \ideal$ for all $j \neq i$, and $\ideal = \talg_{\be}(\fie, \B)$.

Suppose $\be = -\tfrac{1}{n-3}$. Because $e_{i} \in \ideal$, the $\be_{+} = -2/(n-3)$ eigenspace of $L(e_{i})$ is contained in $\ideal$, so $2e_{i} + (n-1)(e_{j} + e_{i\join j})$, and hence also $e_{j} + e_{i\join j}$, are contained in $\ideal$ for all $j \neq i$. Consequently,
\begin{align}
\begin{split}
(n-3)e_{j} & = ((n-4)e_{j} - e_{i\join j}) + (e_{j} + e_{i\join j}) = (n-3)e_{j} \mlt (e_{j} + e_{i\join j}) + (e_{j} + e_{i\join j}) \in \ideal,
\end{split}
\end{align}
so $e_{j} \in \ideal$ for all $j \neq i$, and $\ideal = \talg_{\be}(\fie, \B)$.
\end{proof}

\begin{remark}
By Lemma \ref{eigenlemma}, the two excluded cases, $\be = 1$ and $\be = -\tfrac{n-1}{n-3}$, in Theorem \ref{stsalgebrasimpletheorem} are those values of $\be$ for which $1$ has multiplicity greater than $1$ as an eigenvalue of $L(e_{i})$.

Example \ref{notsimpleexample} shows that the algebra $(\talg_{1}(\fie, AG(2, 3)), \mlt)$ is not simple. 

Lemma \ref{aglemma} shows that for $n = 9$ and $\be = -\tfrac{n-1}{n-3} = - 4/3$, $(\talg_{\be}(\fie, AG(2, 3)), \mlt)$ is simple.
\end{remark}

\section{Axial structure for Hall triple systems}\label{hallsection}
The results of this section show that when $\B$ is a Hall triple system, then $\talg_{\be}(\fie, \B)$ has the structure of an axial algebra for a $\ztwo$-graded fusion law. This result and \cite[Theorem $6.4$]{Hall-Rehren-Shpectorov} are similar in character and their proofs are based on similar computations. These algebras necessarily have large automorphism groups.
For terminology and definitions regarding axial and decomposition algebras the exposition follows principally \cite{Khasraw-McInroy-Shpectorov} and \cite{DeMedts-Peacock-Shpectorov-VanCouwenberghe}.

A \emph{Hall triple system} on $\setn$ is a Steiner triple system on $\setn$ whose underlying Steiner quasigroup is \emph{distributive}, meaning 
\begin{align}\label{hallid1}
&(i\join j)\join (i \join k) = i\join (j \join k),& &\text{pairwise distinct}\,\, i, j, k \in \setn.
\end{align}
Note that \eqref{hallid2} is equivalent to
\begin{align}\label{hallid2}
&k\join((k\join j)\join i) = (k\join i)\join j,& &\text{pairwise distinct}\,\, i, j, k \in \setn. 
\end{align}
The set of Hall triple systems on $\setn$ is denoted $HTS(n)$. Equivalently, a Hall triple system is a Steiner triple system $\B$ for which any two distinct incident blocks generate a subsystem isomorphic to the affine plane of order $3$. Considering the blocks \eqref{9124} of the affine place of order $3$ shows that $\B$ be a Hall triple system means that if $\{i, j, i\join j\}$ and $\{i, k, i \join k\}$ are blocks of $\B$ then so two are the triples lying on every row, column, and left and right diagonal of the tiling of the plane obtained by repeating the diagram in Figure \ref{hallblocks}.
\begin{figure}
\begin{tabular}{c|c|c|c|c|c|c}
  & $\vdots$ & $i \join k$ & $k\join(i \join j)$ & $j\join k$& $\vdots$ &\\
\hline
$\dots$ & $i \join j$& $i$ & $j$ & $i \join j$ & $i$ & $\dots$\\
\hline
$\dots$ & $j \join (i \join k)$ & $k$ & $i \join (j\join k)$ & $j \join (i \join k)$& $k$ & $\dots$\\
\hline
$\dots$  & $j\join k$ & $i \join k$ & $k\join(i \join j)$ & $j\join k$& $i\join k$ & $\dots$\\
\hline
  & $\vdots$ & $i$ & $j$ & $i \join j$ & $\vdots$ &\\
\end{tabular}
\caption{Rows, columns, and diagonals yield blocks}\label{hallblocks}
\end{figure}
See \cite[Chapter $28$]{Colbourn-Dinitz} for a summary of results about Hall triple systems and references.
A $HTS(n)$ is nonempty if and only if $n= 3^{m}$ for some $m \geq 2$ (see \cite{Hall-order} for an elementary proof and references). Because the $m$-dimensional affine space over the field of $3$ elements is an HTS of order $3^{m}$, the content is the claim that the order of a Hall triple system is necessarily a power of $3$.
It is know that up to isomorphism the unique Hall triple system of order $9$ is the affine plane of order $3$, and, more generally, that the smallest order for which there exists a Hall triple system that is not affine space is $81$ \cite{Hall-steiner}. A Steiner triple system with doubly transitive automorphism group is an affine space over the finite field with $3$ elements or a projective space over the finite field with $2$ elements \cite{Hall-steiner}, a Hall triple system is isomorphic to $AG(m, 3)$ for some $m$ if and only if its automorphism group is doubly transitive \cite[Theorem $1$]{Kantor-homogeneous}, and a Steiner triple system is a Fischer space if and only if each triple of points not contained in a block generates an affine plane of order $3$ (see \cite[Section $18$]{Aschbacher}). It is known that a nonaffine Hall triple system of order $3^{m}$ exists for any $m \geq 4$ \cite[Example $28.2$]{Colbourn-Dinitz}, but it seems that no estimate of the number of nonisomorphic Hall triple systems of a given order is known. For combinatorial characterizations of Hall triple systems see \cite{Kral-Macajova-Por-Sereni}. 

A \emph{Fischer space} is a partial linear space system such that the subspace generated by a pair of distinct incident blocks is isomorphic to the dual affine plane of order $2$ or the affine plane of order $3$ (see \cite[Section $18$]{Aschbacher} for the definitions of these spaces). The set of $3$-transpositions of a $3$-transposition group $G$ equipped with the incidence relation comprising blocks of the form $\{x, y, xyx\}$ generates a Fischer space that is a partial Steiner triple system if and only if $G$ is centerless. The Hall triple systems correspond to the subset of Fischer spaces associated with centerless $3$-transposition groups \cite{Cuypers-Hall}. 

For a finite set $\fusion$ with power set $2^{\fusion}$, a \emph{fusion law} over $\fie$ is a symmetric map $\Phi:\fusion \times \fusion \to 2^{\fusion}$. For $n \geq 3$ and $\be \in \fie \setminus\{-\tfrac{n-1}{n-3}, 0, 1\}$ define a fusion law $(\fusion_{\be}, \Phi_{\be})$ by $\fusion_{\be} = \{1, \be_{+}, \be_{-}\}$ and $\Phi_{\be}$ defined as in Table \ref{fusiontable}.
\begin{table}
\begin{tabular}{c||c|c|c}
$\Phi_{\be}$& $1$ & $\be_{+}$ & $\be_{-}$\\
\hline\hline
$ 1$ & $\{1\}$ & $\{\be_{+}\}$ & $\{\be_{-}\}$ \\
\hline
$\be_{+}$  & $\{\be_{+}\}$ & $\{1, \be_{+}\}$ & $\{\be_{-}\}$  \\
\hline
$\be_{-}$ &$\{\be_{-}\}$ & $\{\be_{-}\}$  & $\{1, \be_{+}\}$\\
\end{tabular}
\caption{Fusion law $(\fusion_{\be}, \Phi_{\be})$ for $\be\in \fie\setminus\{-\tfrac{n-1}{n-3}, 0, 1\}$.}\label{fusiontable}
\begin{tabular}{c||c|c|c}
$\Phi_{1/(n-1)}$& $1$ & $0$ & $-2/(n-1)$\\
\hline\hline
$ 1$ & $\{1\}$ & $\emptyset$ & $\{-2/(n-1)\}$ \\
\hline
$0$  & $\emptyset$ & $\{0\}$ & $\{-2/(n-1)\}$  \\
\hline
$-2/(n-1)$ &$\{-2/(n-1)\}$ & $\{-2/(n-1)\}$  & $\{1, 0\}$\\
\end{tabular}
\caption{Fusion law $(\fusion_{1/(n-1)}, \Phi_{1/(n-1)})$ for $\be = \tfrac{1}{n-1}$.}\label{jordanfusiontable}
\end{table}

A morphism $\phi$ of fusion laws $(\fusion_{1}, \Phi_{1}) \to (\fusion_{2}, \Phi_{2})$ is a map of sets $\phi:\fusion_{1}\to \fusion_{2}$ such that $\phi(\fusion_{1}(x, y)) \subset \Phi_{2}(\phi(x), \phi(y))$.
Any finite group $G$ admits the \emph{group fusion law} $(G, \Phi_{G})$ defined by $\Phi_{G}(g, h) = \{gh\}$. A \emph{$G$-grading} of a fusion law $(\fusion, \Phi)$ is a morphism of fusion laws $\phi:(\fusion, \Phi) \to (\fusion_{G}, \Phi_{G})$. 

The fusion law $(\fusion_{\be}, \Phi_{\be})$ of Table \ref{fusiontable} is $\ztwo$-graded if $\be\in \fie\setminus\{-\tfrac{n-1}{n-3}, 0, 1\}$. The morphism $\phi:(\fusion_{\be}, \Phi_{\be}) \to (\ztwo, \Phi_{\ztwo})$ is defined by $\phi(1) = 1$, $\phi(\be_{+}) = 1$, and $\phi(\be_{-}) = -1$ where $\ztwo = \{\pm 1\}$ as a set. The three excluded values of $\be$ are those for which $1$, $\be_{+}$, and $\be_{-}$ are not pairwise distinct.

Given a fusion law $(\fusion, \Phi)$ and a commutative algebra $(\balg, \mlt)$, an element $e \in \balg$ is a \emph{$(\fusion, \Phi)$-axis} if $e$ is an idempotent, $L(e)$ is semisimple with eigenvalues contained in $\fusion$, and the $L(e)$ eigenspaces $\balg^{(\la)}$ satisfy $\balg^{(\la)}\mlt\balg^{(\mu)} \subset \oplus_{\ga \in \Phi(\la, \mu)} \balg^{(\ga)}$. A $(\fusion, \Phi)$-axis $e$ is \emph{primitive} if $1$ has multiplicity $1$ as an eigenvalue of $L(e)$. A \emph{$(\fusion, \Phi)$-axial algebra} is a commutative algebra $(\balg, \mlt)$ equipped with a set $\axes$ of $(\fusion, \Phi)$-axes that generate $(\balg, \mlt)$ as an algebra. A $(\fusion, \Phi)$-axial algebra $((\balg, \mlt), \axes)$ is \emph{primitive} if each element of $\axes$ is primitive.

\begin{theorem}\label{htstheorem}
Consider a Hall triple system $\B \in HTS(n)$.  
\begin{enumerate}
\item For $\be \in \fie\setminus\{-\tfrac{n-1}{n-3}, 0, 1\}$ and the fusion law $(\fusion_{\be}, \Phi_{\be})$ of Table \ref{fusiontable}, the set $\axes = \{e_{i} \in \talg_{\be}: i \in \setn\}$ comprises primitive $(\fusion_{\be}, \Phi_{\be})$-axes in $(\talg_{\be}, \mlt)$, and $((\talg_{\be}(\fie, \B), \mlt), \axes)$ is a primitive $(\fusion_{\be}, \Phi_{\be})$-axial algebra.
\item For $\be = \tfrac{1}{n-1}$ (so $\be_{+} = 0$) and the Jordan fusion law $(\fusion_{1/(n-1)}, \Phi_{1/(n-1)})$ of Table \ref{jordanfusiontable}, the set $\axes = \{e_{i} \in \talg_{\be}: i \in \setn\}$ comprises primitive $(\fusion_{\be}, \Phi_{\be})$-axes in $(\talg_{\be}, \mlt)$, and $((\talg_{\be}(\fie, \B), \mlt), \axes)$ is a primitive $(\fusion_{\be}, \Phi_{\be})$-axial algebra.
\end{enumerate}
\end{theorem}

\begin{proof}
Let $\balg^{e_{i}}_{-}$ and $\balg^{e_{i}}_{+}$ be the $\be_{-}$ and $\be_{+}$ eigenspaces of $L(e_{i})$, as in Lemma \ref{eigenlemma}.
Let $i, j, k, l\in \setn$ be pairwise distinct. Define
\begin{align}\label{fusedeigen}
\balg^{e_{i}}_{+, 1} = \balg^{e_{i}}_{+}+ \fie e_{i} = \fie\{e_{j} + e_{j\join i}: i \neq j \in \setn\} .
\end{align}
Write $\al = (\be -1)/(n-2)$. Straightforward computations using \eqref{hallid1} and \eqref{hallid2} show
\begin{align}\label{axial}
\begin{aligned}
(e_{j} &- e_{j\join i})\mlt(e_{j} - e_{j\join i}) = (1 - 2\al)(e_{j} + e_{j\join i}) - 2\be e_{i} \in \balg^{e_{i}}_{+, 1} ,\\
(e_{j} &+ e_{j\join i})\mlt(e_{j} + e_{j\join i}) = (1 + 2\al)(e_{j} + e_{j\join i}) +  2\be e_{i}\in \balg^{e_{i}}_{+, 1},\\
(e_{j} &- e_{j\join i})\mlt(e_{j} + e_{j\join i}) = e_{j} - e_{j\join i}\in \balg^{e_{i}}_{-} ,
\end{aligned}
\end{align}
\begin{align}\label{axial2}
\begin{aligned}
(e_{j} &- e_{j\join i})\mlt(e_{k} - e_{k\join i}) \\
& = \be(e_{j\join k} + e_{(i\join j)\join (i \join k)} - e_{(i\join j)\join k} - e_{(i \join k)\join j})\\
& = \be(e_{j\join k} + e_{(j\join k)\join i} - e_{(i\join k)\join j} - e_{((i\join k)\join j)\join i}) \in \balg^{e_{i}}_{+},\\
(e_{j} &+ e_{j\join i})\mlt(e_{k} + e_{k\join i}) \\
&=2\al(e_{j} + e_{j\join i} + e_{k} + e_{k\join i}) + \be(e_{j\join k} + e_{(i\join j)\join (i \join k)} + e_{(i\join j)\join k} + e_{(i \join k)\join j})\\
& = 2\al(e_{j} + e_{j\join i} + e_{k} + e_{k\join i}) + \be (e_{j\join k}  + e_{(j\join k)\join i} + e_{(i\join k)\join j}  + e_{((i\join k)\join j)\join i}) \in \balg^{e_{i}}_{+, 1},\\
(e_{j} &- e_{j\join i})\mlt(e_{k} + e_{k\join i}) \\
&= 2\al(e_{j} - e_{j\join i}) + \be(e_{j\join k} - e_{(i\join j)\join (i \join k)} - e_{(i\join j)\join k} + e_{(i \join k)\join j})\\
& =2\al(e_{j} - e_{j\join i}) + \be(e_{j\join k}  - e_{(j\join k)\join i} + e_{(i\join k)\join j} - e_{((i\join k)\join j)\join i}) \in \balg^{e_{i}}_{-}.
\end{aligned}
\end{align}
Because $(e_{j} \pm e_{j\join i}) \mlt e_{i} =\be_{\pm}(e_{j} - e_{j\join i}) \in \balg^{e_{i}}_{\pm}$, together with \eqref{alplusbeeigen}, \eqref{alminusbeeigen} and \eqref{fusedeigen}, the relations \eqref{axial} suffice to show that $\balg^{e_{i}}_{-}\mlt \balg^{e_{i}}_{-} \subset\balg^{e_{i}}_{+, 1}$, $\balg^{e_{i}}_{-}\mlt \balg^{e_{i}}_{+} \subset \balg^{e_{i}}_{-}$, and $\balg^{e_{i}}_{+}\mlt \balg^{e_{i}}_{+} \subset \balg^{e_{i}}_{+, 1}$. Together with Lemma \ref{eigenlemma} this shows that for $\be \in \fie \setminus\{-\tfrac{n-1}{n-3}, 0, 1\}$, $e_{i}$ is a primitive $(\fusion_{\be}, \Phi_{\be})$-axis. Because by definition $\axes$ generates $(\talg_{\be}(\fie, \B), \mlt)$, this shows that $((\talg_{\be}(\fie, \B), \mlt), \axes)$ is a primitive $(\fusion_{\be}, \Phi_{\be})$-axial algebra.

There remains to consider the particular case $\be = \tfrac{1}{n-1}$, in which case $\be_{+} = 0$. For general $\be$, computations using \eqref{axial} and \eqref{axial2} show
\begin{align}\label{axial3}
\begin{aligned}
(2e_{i} &+ (n-1)(e_{j} + e_{i \join j}))\mlt (2e_{i} + (n-1)(e_{j} + e_{i \join j})) \\
&= (n-3)(2e_{i} + (n-1)(e_{j} + e_{i \join j})) + \be_{+}((n^{2} - 3n + 6)e_{i} + 6(n-1)(e_{j} + e_{i\join j}),\\
(2e_{i} &+ (n-1)(e_{j} + e_{i \join j}))\mlt (2e_{i} + (n-1)(e_{k} + e_{i \join k})) \\
&= - 4n\be_{+}e_{i} + 2(2\be_{+} - 1)(2e_{i} + (n-1)(e_{j} + e_{i \join j}) + 2e_{i} + (n-1)(e_{k} + e_{i \join k}))\\
& + (n-1)\be\left(2e_{i} + (n-1)(e_{j\join k}  + e_{(j\join k)\join i}) + 2e_{i} + (n-1)(e_{(i\join k)\join j}  + e_{((i\join k)\join j)\join i}) \right).
\end{aligned}
\end{align}
Together with \eqref{alplusbeeigen} of Lemma \ref{eigenlemma}, the relations \eqref{axial3} show $\balg^{e_{i}}_{+}$ is a subalgebra if and only if $\be_{+} = 0$. Because $e_{i}\mlt(2e_{i} + (n-1)(e_{j} + e_{i \join j})) = \be_{+}(2e_{i} + (n-1)(e_{j} + e_{i \join j}))$, it follows that, when $\be_{+} = 0$, $((\talg_{\be}(\fie, \B), \mlt), \axes)$ is a primitive $(\fusion_{\be}, \Phi_{\be})$-axial algebra for the Jordan fusion rule of Table \ref{jordanfusiontable}. Moreover, it follows the axes of $((\talg_{\be}(\fie, \B), \mlt), \axes)$ do not satisfy a Jordan fusion rule when $\be_{+} \neq 0$.
\end{proof}

For $\be \in \fie \setminus\{-\tfrac{n-1}{n-3}, 0, 1\}$ and $\B \in STS(n)$ the endormophism $\miya_{i} = \miya_{e_{i}}\in \eno(\talg_{\be}(\fie, \B))$ defined by extending linearly the map satisfying $\miya_{i}(e_{i}) = e_{i}$ and $\miya_{i}(e_{j}) = e_{i\join j}$ for $i \neq j \in \setn$ satisfies $\miya_{i}(e_{j}\pm e_{j\join i}) = \pm (e_{j}\join e_{j\join i})$. This shows $\miya_{i}(\balg^{e_{i}}_{+, 1}) = \balg^{e_{i}}_{+, 1}$ and $\miya_{i}(\balg^{e_{i}}_{-}) = \balg^{e_{i}}_{-}$ so that $\miya_{i}$ is the reflection through $\balg^{e_{i}}_{+, 1}$ along $\balg^{e_{i}}_{-}$. It is straightforward to check that $\miya_{i}$ is an automorphism of $(\talg_{\be}(\fie, \B), \mlt)$ if and only if there hold the equivalent identities \eqref{hallid1} and \eqref{hallid2}, that is if and only if $\B$ is a Hall triple system. In the case $\B \in HTS(n)$, it follows from Theorem \ref{htstheorem} that $\miya_{i}$ is the automorphism of $(\talg_{\be}(\fie, \B), \mlt)$ called the \emph{Miyamoto involution} determine by the axis $e_{i}$ and the $\ztwo$-grading of the fusion law $(\fusion_{\be}, \Phi_{\be})$ \cite{Khasraw-McInroy-Shpectorov}.

The group $G$ generated by $\{\miya_{i}: i \in \setn\}$ is the \emph{Miyamoto group} of $((\talg_{\be}(\fie, \B), \mlt), \axes)$. Because $\miya_{i}^{2} = \Id$ and $\miya_{i}\circ \miya_{j}\circ \miya_{i} = \miya_{i\join j} = \miya_{j}\circ \miya_{i}\circ \miya_{j}$, the set $\{\miya_{i}:i \in \setn\}$ is a single conjugacy class in $G$ comprising involutions such that any product $\miya_{i}\circ \miya_{j}$ has order three \cite{Hall-automorphisms}. That is $\{\miya_{i}: i \in \setn\}$ is a conjugacy class of $3$-transpositions \cite{Aschbacher, Cuypers-Hall} and $G$ is a \emph{Fischer group} \cite[Definiton $3.2$]{Manin-cubicforms}. Because $\miya_{i}\miya_{j} =  \miya_{i\join j}\circ \miya_{j}\circ \miya_{i\join j}\circ \miya_{j}$ is a commutator, the order of the abelianization $G/[G, G]$ is at most two. By theorems of Fischer (see \cite{Fischer} and \cite[Theorems $9.1$ and $9.2$]{Manin-cubicforms}), $G/[G, G] \simeq \ztwo$ and $[G, G]$ is a nilpotent $3$-group.

Via the identification of the set of axes $\axes$ with $\setn$, each $\miya_{i}$ acts as the involutory automorphism $\si_{i} \in S_{n}$ of $\B$ defined by $\si_{i}(j) = i \join j$ for $j \in \setn$. The Miyamoto group $G$ is evidently isomorphic to the group $K \subset S_{n}$ generated by these automorphisms $\si_{i}$ of $\B$, and, by \cite[Lemma $4.1$]{Hall-automorphisms}, $K$ acts transitively on $\B$ and $\si_{i}$ centralizes the stabilizer of $i$ in $K$. This shows that the Miyamoto group $G$ acts as automorphisms of $\B$ and hence also of $(\talg_{\be}(\fie, \B), \mlt)$. Because the Miyamoto group of $(\talg_{\be}(\fie, \B), \mlt)$ is nontrivial, the automorphism group of $(\talg_{\be}(\fie, \B), \mlt)$ is never trivial. 

In general the Miyamoto group is a proper subgroup of $\Aut(\B)$. For example for $\B = AG(2, 3)$, the Miyamoto group $G$ has order $54$ (it is generated by any $\tau_{i}$, $\tau_{j}$, and $\tau_{k}$ such that $i \neq j$ and $k \notin \{i,j, i\join j\}$; see \cite[p. $469$]{Hall-automorphisms}), while $\Aut(\B)$ is the affine group of a two-dimensional affine space over the finite field with $3$ elements, which has order $432$.

Whether any two of the Miyamoto group, $\Aut(\B)$, and $\Aut(\talg_{\be}(\fie, \B), \mlt)$ are isomorphic has not been established in general. Theorem \ref{automorphismtheorem} shows that in a special case the automorphism groups of $\B$ and $(\talg_{\be}(\fie, \B), \mlt)$ are isomorphic. Its proof uses Lemma \ref{permutelemma}. The argument depends strongly on the assumption that the base field is Euclidean and works only for certain values of $\be$.

\begin{lemma}\label{permutelemma}
Let $\rea$ be a Euclidean field and consider a Hall triple system $\B \in HTS(n)$. Let $\be \in \rea$ satisfy $-\tfrac{1}{n-3} < \be < \tfrac{1}{n-1}$. If $f \in \talg_{\be}(\rea, \B)$ is idempotent and satisfies $\kform(f, f) = \kform(e_{j}, e_{j})$ for any $j \in \setn$, then $f = e_{i}$ for some $i \in \setn$.
\end{lemma}

\begin{proof}
When $\be = 0$, $(\talg_{0}(\rea, \B), \mlt) \simeq \ealg^{n-1}(\rea)$ and the claim follows from the complete description of the idempotents in $\ealg^{n-1}(\rea)$ due to \cite[Lemma $2$ and Corollary $3$]{Harada} (see \cite[Lemma $3.16$]{Fox-simplicial}; the proofs are given for the real field, but work without change over a general Euclidean field) which states that every idempotent is a multiple of a sum of the form $\sum_{i \in I}e_{i}$ for some $I \subset \setn$ and computes $\kform$ on these elements. For the rest of the proof, suppose $\be \neq 0$.

For $i \in \setn$, by Theorem \ref{htstheorem}, $f = \la e_{i} + f_{+} + f_{-}$ where $f_{\pm} \in \balg^{e_{i}}_{\pm}$. By the invariance of $\kform$, the subspaces $\fie e_{i}$, $\balg^{e_{i}}_{+}$, and $\balg^{e_{i}}_{-}$ are $\kform$-orthogonal, so, from $\kform(e_{i}, e_{i}) = \kform(f, f)$, there follows
\begin{align}\label{perm0}
(1-\la^{2})\kform(e_{i}, e_{i}) =  \kform(f_{+}, f_{+})+\kform(f_{-}, f_{-}).
\end{align} 
By Theorem \ref{htstheorem}, decomposing $f\mlt f = f$ into its projections onto $\balg^{e_{-}}_{1, +}$ and $\balg^{e_{i}}_{-}$ yields
\begin{align}\label{perm1}
&\la(1-\la) e_{i} =  (2\la \be_{+}-1)f_{+} + f_{+}\mlt f_{+} + f_{-}\mlt f_{-},& & (1-2\la\be_{-})f_{-} = 2f_{+}\mlt f_{-}.
\end{align}
Pairing the first equation of \eqref{perm1} with $e_{i}$ and using the invariance of $\kform$ yields 
\begin{align}\label{perm2}
\begin{split}
\la(1-\la)\kform(e_{i}, e_{i}) &= \kform(e_{i},  f_{+}\mlt f_{+} ) + \kform(e_{i}, f_{-}\mlt f_{-}) 
= \be_{+}\kform(f_{+}, f_{+})+ \be_{-}\kform(f_{-}, f_{-}).
\end{split}
\end{align}
By Theorem \ref{stsalgebratheorem}, $\kform$ is positive definite, so \eqref{perm0} implies $\la \in (-1, 1]$ (were $\la = -1$, then both $f$ and $-f$ would be idempotent).
The assumption $-\tfrac{1}{n-3} < \be < \tfrac{1}{n-1}$ is equivalent to $\be_{+}$ and $\be_{-}$ both being negative. By \eqref{perm2}, this implies $\la(1-\la)\kform(e_{i}, e_{i}) \leq 0$, which, because $\la \leq 1$, implies either $\la = 1$ or $\la \leq 0$. If $\la = 1$, then $f = e_{i}$ and there is nothing more to show, so suppose $\la \leq 0$, in which case $\la \in (-1, 0]$. 
Let $i \in \setn$ satisfy $\kform(f, e_{i}) = \max\{\kform(f, e_{j}): j \in \setn\}$. 
Were $\kform(f, e_{i}) < 0$, then $\sum_{j \neq i}\kform(f, e_{j})  =- \kform(f, e_{i}) > 0$, which would imply there is $j \in \setn$ such that $\kform(f, e_{j}) > 0 > \kform(f, e_{i})$ contrary to the choice of $i$. Hence $\la \kform(e_{i}, e_{i}) = \kform(f, e_{i}) \geq 0$ implies $\la \geq 0$. Since $\la \in (-1, 0]$, $\la = 0$. In \eqref{perm2} this yields $f_{+} =0$ and $f_{-} = 0$, so $f = 0$, a contradiction with $\kform(f, f) = \kform(e_{i}, e_{i}) \neq 0$.
\end{proof}

\begin{theorem}\label{automorphismtheorem}
For a Hall triple system $\B \in HTS(n)$ and $\be \in (-\tfrac{1}{n-3}, 0) \cup (0, \tfrac{1}{n-1}) \subset \rea$, $\Aut(\talg_{\be}(\rea, \B), \mlt)$ is isomophic to $\Aut(\B)$. 
\end{theorem}

\begin{proof}
Because $\phi \in \Aut(\talg_{\be}(\rea, \B), \mlt)$ preserves $\kform$, by Lemma \ref{permutelemma}, $\phi$ preserves $\{e_{j}: j \in \setn\}$ so there is a permutation $\si \in S_{n}$ such that $\phi(e_{i}) = e_{\si(i)}$ for $i \in \setn$. By Lemma \ref{prepermutelemma}, $\si$ is an automorphism of $\B$. This shows that the canonical group embedding $\Aut(\B) \to \Aut(\talg_{\be}(\rea, \B), \mlt)$ is surjective.
\end{proof}

The value $\be = -\tfrac{n-1}{n-3}$ in Lemma \ref{htsalgebrasimplelemma} is one of the values excluded in Theorem \ref{stsalgebrasimpletheorem}.

\begin{lemma}\label{htsalgebrasimplelemma}
Consider a Hall triple system $\B \in HTS(n)$. For $\be =  -\tfrac{n-1}{n-3}$:
\begin{enumerate}
\item The algebra $(\talg_{\be}(\fie, \B), \mlt)$ contains no proper, nontrivial ideal stable under any Miyamoto involution.
\item\label{htssim2} If the algebra $(\talg_{\be}(\fie, \B), \mlt)$ contains a proper, nontrivial ideal $\ideal$, then $\dim \ideal = (n-1)/2$, the $\kform$-orthogonal complement $\ideal^{\perp}$ is an ideal transverse to $\ideal$, and $\ideal$ and $\ideal^{\perp}$ are each stable under the action of the commutator subgroup of the Miyamoto group, while any Miyamoto involution interchanges $\ideal$ and $\ideal^{\perp}$.
\end{enumerate}
\end{lemma}

\begin{proof}
For $i \in \setn$, let $\miya_{i}$ be the Miyamoto involution associated with $e_{i}$. Because $\miya_{i}$ is an automorphism, if $\ideal \subset \talg_{\be}(\fie, \B)$ is a nontrivial ideal, $\miya_{i}(\ideal)$ is also a nontrivial ideal. Let $\ideal \subset \talg_{\be}(\fie, \B)$ be a proper, nontrivial ideal satisfying $\miya_{i}(\ideal) = \ideal$ and let $x = \sum_{j \in \setn}x_{j}e_{j} \in \ideal$ be nonzero. Combining
\begin{align}
&L(e_{i})x - \al x = ((1-\al)x_{i} + \al \sum_{j \neq i}x_{j})e_{i} + \be\sum_{j \neq i}x_{i\join j}e_{j}, & &\miya_{i}(x) = x_{i}e_{i} + \sum_{j \neq i}x_{i\join j}e_{j},&
\end{align}
yields
\begin{align}\label{htss1}
 ((1-\be_{+})x_{i} + \al \sum_{j \neq i}x_{j})e_{i} & = L(e_{i})x - \al x - \be \miya_{i}(x)  \in \ideal.
 \end{align}
 If $\be_{+} \neq 0$, the left-hand side of \eqref{htss1} vanishes if and only if all the $x_{j}$ are equal, in which case $x = 0$, contrary to hypothesis. Consequently, if $\be_{+} \neq 1$, then $e_{i} \in \ideal$. When $\be = -\tfrac{n-1}{n-3}$, $L(e_{i})$ is invertible by Lemma \ref{eigenlemma}, so $e_{i} \in \ideal$ implies $\ideal = \talg_{\be}(\fie, \B)$. 
 
Let $\ideal \subset \talg_{\be}(\fie, \B)$ be a nontrivial ideal. Were the $\miya_{i}$-stable ideal $\ideal \cap \miya_{i}(\ideal)$ nontrivial, then, by the preceding, it, and hence also $\ideal$ would equal $\talg_{\be}(\fie, \B)$. Hence if $\ideal \subset \talg_{\be}(\fie, \B)$ is a proper, nontrivial ideal then $\ideal \cap \miya_{i}(\ideal) = \{0\}$. In this case, again by the preceding, the nontrivial $\miya_{i}$-stable ideal $\ideal + \miya_{i}(\ideal)$ must equal $\talg_{\be}(\fie, \B)$. This shows that if $\ideal$ is a proper, nontrivial ideal, then $\talg_{\be}(\fie, \B) = \ideal \oplus \miya_{i}(\ideal)$ for all $i \in \setn$.  If $a, b \in \ideal$, $L(a)$ preserves $\ideal$ and annihilates $\miya_{i}(\ideal)$ while $L(\miya_{i}(b))$ preserves $\miya_{i}(\ideal)$ and annihilates $\ideal$. It follows that $L(a)L(\miya_{i}(b))$ is identically $0$. In particular $\kform(a, \miya_{i}(b)) = 0$ and the direct sum $\talg_{\be}(\fie, \B) = \ideal \oplus \miya_{i}(\ideal)$ is $\kform$-orthogonal. Thus $\miya_{i}(\ideal)$ is the $\kform$-orthogonal complement of $\ideal$ for all $i \in \setn$, and, consequently, $\miya_{i}(\ideal) = \miya_{j}(\ideal)$ for all $i, j \in \setn$. Equivalently the order three elements $\miya_{ij} = \miya_{i}\circ \miya_{j}$ of the Miyamoto group $G$ stabilize $\ideal$, so the commutator subgroup $[G, G]$ stabilizes $\ideal$. 
\end{proof}

For $\be =  -\tfrac{n-1}{n-3}$, it is not clear if $(\talg_{\be}(\fie, \B), \mlt)$ must be simple for $\B \in HTS(n)$, although Lemma \ref{aglemma} shows that for $\be = -4/3$, $(\talg_{\be}(\fie, AG(2, 3)), \mlt)$ is simple, and the author expects that the same is true at least for the Hall triple systems determined by the affine spaces $AG(m, 3)$. It would suffice to show that an ideal as in \eqref{htssim2} of Lemma \ref{htsalgebrasimplelemma} cannot exist. It would be interesting to know also if, for $\be =  -\tfrac{n-1}{n-3}$, whether $(\talg_{\be}(\fie, \B), \mlt)$ can fail to be simple for a non-Hall $\B \in STS(n)$.

\section{Partial results about idempotents and square-zero elements}\label{idempotentsection}
The remainder of the paper records some partial results about the structure of idempotents and square-zero elements in $\talg_{\be}(\fie, \B)$.

A block $B \in \B$ containing $i \in \setn$ can be written $B = \{i, j, i \join j\}$. Although the choice of $j$ is ambiguous, the associated quadratic polynomial 
\begin{align}
Q_{i, B}(x) = (\be_{+} + \be_{-})x_{i}(x_{j} + x_{i\join j}) + (\be_{+} - \be_{-})x_{j}x_{i\join j} = \tfrac{2(\be-1)}{n-2}x_{i}(x_{j} + x_{i\join j})+ 2\be x_{j}x_{i\join j} ,
\end{align}
depends only on the pair $(i, B)$. 

\begin{lemma}
For $\ep \in \{0, 1\}$, $\be \in \fie$, and $\B \in STS(n)$, an element $x =\sum_{i \in \setn}x_{i}e_{i}\in \talg_{\be}(\fie, \B)$ satisfies $x \mlt x = \ep x$ if and only if there is $c \in \fie$ such that for all $i \in \setn$,
\begin{align}\label{talgidemequations}
\begin{split}
x_{i}^{2} - \ep x_{i} + \sum_{i \in B}Q_{i, B}(x) = c.
\end{split}
\end{align}
In this case, writing $\bar{x} = \tfrac{1}{n}\sum_{j \in \setn}x_{j}$,
\begin{align}\label{nc}
\begin{split}
c &
= \tfrac{1-\be_{+}}{n-1}\sum_{j \in \setn}(x_{j} - \bar{x})^{2} - \ep \bar{x} + n\be_{+}\bar{x}^{2}= \tfrac{1-\be}{n-2}\sum_{j \in \setn}(x_{j} - \bar{x})^{2} - \ep \bar{x} +\tfrac{n((n-1)\be - 1)}{n-2}\bar{x}^{2} .
\end{split}
\end{align}
\end{lemma}
\begin{proof}
Note that the coefficients $x_{i}$ are determined only up to adding simultaneously to all of them some fixed constant. However, the symmetry of the representation with respect to the full spanning set $\{e_{i}: i \in \setn\}$ is more convenient than the representation of $x$ with respect to a basis contained in this spanning set. For $x \in \talg_{\be}$,
\begin{align}\label{talgxx}
\begin{split}
x\mlt x & = \sum_{i \in \setn}x_{i}^{2}e_{i} +  \sum_{i \in \setn}\sum_{j \neq i}x_{i}x_{j}\left(\tfrac{\be - 1}{n-2}(e_{i} +e_{j}) + \be  e _{i\join j}\right)\\
&  = \sum_{i \in \setn}\left( x_{i}^{2} +\tfrac{2(\be - 1)}{n-2}\sum_{j \neq i}x_{i}x_{j} + \be \sum_{j \neq i} x_{j}x_{i \join j}\right)e_{i} \\
&= \sum_{i \in \setn}\left( x_{i}^{2} +(\be_{+} + \be_{-})\sum_{j \neq i}x_{i}x_{j} + \tfrac{\be_{+} - \be_{-}}{2}\sum_{j \neq i} x_{i\join j}x_{j} \right)e_{i} 
= \sum_{i \in \setn}\left(x_{i}^{2} + \sum_{i \in B}Q_{i, B}(x) \right)e_{i}.
\end{split}
\end{align}
By \eqref{talgxx}, $x \mlt x = \ep x$ if and only if there is $c \in \rea$ such that \eqref{talgidemequations} holds
for $i \in \setn$. Summing \eqref{talgidemequations} using $\sum_{j \neq i}x_{i\join j} = \sum_{j \neq i}x_{j} = n\bar{x} - x_{i}$, $\sum_{j \in \setn}x_{j}^{2}= \sum_{j \in \setn}(x_{j} - \bar{x})^{2} + n\bar{x}^{2}$, and
\begin{align}
\begin{aligned}
\sum_{i \in \setn}\sum_{j \neq i} x_{i\join j}x_{j} &=\sum_{j \in \setn}\sum_{i \neq j} x_{i\join j}x_{j} =  \sum_{j \in \setn}\sum_{i \neq j} x_{i}x_{j} \\&= n^{2}\bar{x}^{2} - \sum_{j \in \setn}x_{j}^{2} = n(n-1)\bar{x}^{2} - \sum_{j \in \setn}(x_{j} - \bar{x})^{2}, 
\end{aligned}
\end{align}
yields
\begin{align}\label{talgidemequations1}
\begin{split}
nc 
& = \sum_{i \in \setn}x_{i}^{2} - \ep n\bar{x}+ (\be_{+} + \be_{-})\sum_{i \in \setn}\sum_{j \neq i}x_{i}x_{j} + \tfrac{\be_{+} - \be_{-}}{2}\sum_{i \in \setn}\sum_{j \neq i}x_{j}x_{i\join j}\\
& = - \ep n\bar{x}+ (1 - \tfrac{3\be_{+} + \be_{-}}{2})\sum_{j \in \setn}(x_{j} - \bar{x})^{2}  + \tfrac{n}{2} \left(1 + (n-1)(3\be_{+} + \be_{-})\right)\bar{x}^{2}\\
& = -\tfrac{n(\be_{+} + \be_{-})}{2}\sum_{j \in \setn}(x_{j} - \bar{x})^{2} - \ep n \bar{x} + n^{2}\be_{+}\bar{x}^{2},
\end{split}
\end{align}
in which there have been used $(n-1)(3\be_{+} + \be_{-}) = 2n\be_{+} - 2$ and $2 - (3\be_{+} + \be_{-}) = n(\be_{+} + \be_{-})$.
This shows \eqref{nc}. 
\end{proof}

Fix a reference block $B(i, j) = \{i, j, i \join j\} \in \B$. Taking the difference of \eqref{talgidemequations} with itself with $j$ in place of $i$ yields 
\begin{align}
\label{talgidemdiffequations}
\begin{split}
0 &= (x_{i} -x _{j})\left(x_{i} + x_{j}- 2\be x_{i\join j} - \ep+ \tfrac{2(\be - 1)}{n-2}\sum_{k \notin \{i, j\}}x_{k}  \right) + \be\sum_{k \notin B(i, j) }x_{k}(x_{i\join k} - x_{j\join k})\\
&= (x_{i} -x _{j})\left(x_{i} + x_{j}  +2\be_{-}x_{i\join j}  - \ep + \tfrac{2(\be - 1)}{n-2}\sum_{k \notin B(i, j) }x_{k}\right) +\be \sum_{k \notin B(i, j) }x_{k}(x_{i\join k} - x_{j\join k}),
\end{split}
\end{align}
for $1 \leq i \neq j \leq n$.
The equations \eqref{talgidemdiffequations} can be rewritten as
\begin{align}\label{blocktalgidemdiffequations}
\begin{split}
&\resizebox{.9\textwidth}{!}{$\begin{pmatrix*}[c]
x_{j} - x_{i\join j} & 0 & 0\\ 
0 & x_{i\join j} - x_{i} & 0 \\ 
0 & 0 & x_{i} - x_{j}
\end{pmatrix*}
\left(\left((2\be_{-} - 1)\imat_{3} +\ones_{3}\right)
\begin{pmatrix} x_{i}\\x_{j}\\x_{i\join j}\end{pmatrix}
- \left(\ep + \tfrac{2(1-\be)}{n-2}\sum_{k \notin B(i, j)}x_{k}\right)\begin{pmatrix}1\\1\\1\end{pmatrix}
\right)$} \\
&= - \be \sum_{k \notin B(i, j) }x_{k}
\begin{pmatrix*}[c]
x_{j\join k} - x_{(i \join j)\join k}\\
x_{(i\join j)\join k} - x_{i\join k}\\
x_{i\join k} - x_{j\join k}
\end{pmatrix*}
= - \be \sum_{k \notin B(i, j) }x_{k}
\begin{pmatrix*}[r]
0 & 1 & - 1\\
-1 & 0 & 1 \\
1 & 0 & -1 \\\end{pmatrix*}
\begin{pmatrix*}[c]
x_{i\join k} \\ x_{j\join k} \\ x_{(i \join j)\join k}\\
\end{pmatrix*},
\end{split}
\end{align}
where $\imat_{3}$ is the $3 \times 3$ identity matrix and $\ones_{3}$ is the $3\times 3$ matrix with all components equal to $1$. Solving \eqref{blocktalgidemdiffequations} yields Lemma \ref{talgidempotentlemma}, which exhibits idempotents and square-zero elements in $(\talg_{\be}(\fie, \B), \mlt)$ associated with a block $B \in \B$.

For the special case $\B = \{1, 2, 3\}$, the claims of Lemma \ref{talgidempotentlemma} can be found in \cite[Section $5.3$]{Krasnov-Tkachev}.

\begin{lemma}\label{talgidempotentlemma}
Suppose $n > 3$. Let $\be \in \fie$ and $\B \in STS(n)$. For $B \in \B$ let $B^{c} = \{k \in \setn: k \notin B\}$. 
For any $i \neq j \in \setn$ and $B = \{i, j, i\join j\} \in \B$, any idempotent or square-zero element of $(\talg_{\be}(\fie, \B), \mlt)$ contained in $\fie\{e_{i}, e_{j}, e_{i\join j}\}$ is among those appearing in the following list.
\begin{enumerate}
\item\label{ebeb} If $\be \neq \tfrac{6-n}{2n}$, then $e^{0}_{B} = \tfrac{n-2}{2n\be + n-6}(e_{i} + e_{j} + e_{i\join j})$ is idempotent in $(\talg_{\be}(\fie, \B), \mlt)$. 
\item\label{zbzb} If $\be = \tfrac{6-n}{2n}$, then $z_{B} = e_{i} + e_{j} + e_{i\join j}$ is square-zero in $(\talg_{\be}(\fie, \B), \mlt)$.
\item\label{elel} If $\be =  -\tfrac{n}{2(n-3)}$ (in this case $\be_{-} = 1/2$) then $e = \la_{1}e_{i} + \la_{2}e_{j} + \la_{3}e_{i\join j}$ solves $e\mlt e = \ep e$ for $\ep \in \{0, 1\}$ if and only if
\begin{align}\label{elambdab}
&\la_{1} + \la_{2} + \la_{3} = \ep,& &\la_{1}^{2} + \la_{2}^{2} + \la_{3}^{2} = \ep^{2}.
\end{align}
\item If $1 - \be_{+} + \be_{-} - 4\be_{+}\be_{-} = \tfrac{4(n-3)(n-1)\be^{2} - 2n(n-4)\be + n(n-4)}{n-2} \neq 0$, the three elements $e_{B}^{i}$, $e_{B}^{j}$ and $e_{B}^{i\join j}$ defined by
\begin{align}\label{fbidem}
\begin{split}
e_{B}^{i\join j} &  = \tfrac{1}{1 - \be_{+} + \be_{-} - 4\be_{+}\be_{-}}\left((1 - 2\be_{+})(e_{i} + e_{j}) + 2\be e_{i\join j}\right)\\
& = \tfrac{n-2}{4(n-3)(n-1)\be^{2} - 2n(n-4)\be + n(n-4)}\left( (n - 2(n-1)\be)(e_{i} + e_{j}) + 2(n-2)\be e_{i\join j}\right),
\end{split}
\end{align}
and its cyclic permutations in $\{i, j, i\join j\}$ are idempotents in $(\talg_{\be}(\fie, \B), \mlt)$ that are, moreover, linearly independent if $\be \notin \{\tfrac{n}{2(n-1)}, -\tfrac{n}{2(n-3)}\}$.
\begin{enumerate}
\item If $\be = \tfrac{n}{2(n-1)}$, so that $\be_{+} = 1/2$, then $e_{B}^{i \join j} = e_{i\join j}$.
\item If $\be = -\tfrac{n}{2(n-3)}$, so that $\be_{-} = 1/2$, then $e_{B}^{i}  + e_{B}^{j} +e_{B}^{i\join j} = 0$.
\end{enumerate}
\end{enumerate}
\end{lemma}
\begin{proof}
Claims \eqref{ebeb} and \eqref{zbzb} are immediate from \eqref{gabgab}.
The form of \eqref{blocktalgidemdiffequations} suggests looking for $x$ solving $x\mlt x = \ep x$ such that $x_{k} = 0$ for $k \in B^{c}$ (because $\sum_{k \in \setn}e_{k} = 0$, nothing is gained by considering the apparently more general condition that $x_{k} = t$ for $k \in B^{c}$).
Supposing that $x_{k} = 0$ for $k \in B^{c}$ simplifies \eqref{talgidemequations} and \eqref{blocktalgidemdiffequations} considerably to
\begin{align}\label{blocktalgidemequationssimp}
\begin{split}
x_{i}^{2} - \ep x_{i} + (\be_{+} + \be_{-})(x_{i}x_{j} + x_{i}x_{i\join j}) +(\be_{+} - \be_{-})x_{j}x_{i\join j}  = c,
\end{split}
\end{align}
(and those equations obtained from it by cyclically permuting $\{i, j, i \join j\}$)
and
\begin{align}\label{blocktalgidemdiffequationssimp}
\begin{split}
&\begin{pmatrix}
x_{j} - x_{i\join j} & 0 & 0\\ 
0 & x_{i\join j} - x_{i} & 0 \\ 
0 & 0 & x_{i} - x_{j}
\end{pmatrix}
\left(
\left((2\be_{-} - 1)\imat_{3} + \ones_{3}\right)
\begin{pmatrix} x_{i}\\x_{j}\\x_{i\join j}\end{pmatrix}
- \begin{pmatrix}\ep\\\ep\\\ep\end{pmatrix}
\right) = 0
\end{split}
\end{align}
In the rest of the proof, supoose $x_{k} = 0$ for $k \in B^{c}$. If $x_{i} = x_{j} = x_{i\join j} = s \neq 0$, \eqref{blocktalgidemequationssimp} becomes $s((1 +3\be_{+} + \be_{-})s - \ep) = c$. On the other hand, \eqref{nc} becomes $nc = 3s((1 + 3\be_{+} + \be_{-})s - \ep)$. Since $n > 3$, these equations imply $s((1 +3\be_{+} + \be_{-})s - \ep) = 0$. When $\ep = 1$ this admits the nontrivial solution $s = \tfrac{n-2}{2n\be + n-6}$ as long as $\be \neq \tfrac{6-n}{2n}$ (equivalently $1 + 3\be_{+} + \be_{-} \neq 0$), while when $\ep = 0$, it is inconsistent except when $\be = \tfrac{6-n}{2n}$. This recovers the solutions \eqref{ebeb} and \eqref{zbzb}. 

Since all solutions such that $x_{i} = x_{j} = x_{i\join j}$ have been found already, there remains to consider the cases in which $x_{i}$, $x_{j}$, and $x_{i\join j}$ are not all equal.

A solution of
\begin{align}\label{blocktalgidemdiffequationssimp2}
\begin{split}
\left((2\be_{-} - 1)\imat_{3} + \ones_{3}\right)
\begin{pmatrix} x_{i}\\x_{j}\\x_{i\join j}\end{pmatrix}
= \begin{pmatrix}\ep\\\ep\\\ep\end{pmatrix},
\end{split}
\end{align}
yields a solution of \eqref{blocktalgidemdiffequationssimp}. If $x_{i}$, $x_{j}$, and $x_{i\join j}$ are pairwise distinct, $x$ solves \eqref{blocktalgidemdiffequationssimp} if and only if $x$ solves \eqref{blocktalgidemdiffequationssimp2}. 

The eigenvalues of $(2\be_{-} - 1)\imat_{3} + \ones_{3}$ are $2\be_{-} - 1$ and $2(\be_{-} +1)$. If these are nonzero, then \eqref{blocktalgidemdiffequationssimp2} has a unique solution whose components are necessarily equal because $\left((2\be_{-} - 1)\imat_{3} + \ones_{3}\right)^{-1} = \tfrac{1}{2\be_{-} - 1}(\imat_{3} - \tfrac{1}{2(\be_{-} + 1)}\ones_{3})$. Solving \eqref{blocktalgidemdiffequationssimp2} in these cases recovers the solution $e_{B}$ or $z_{B}$ already found and there are no solutions of  \eqref{blocktalgidemdiffequationssimp} with $x_{i}$, $x_{j}$, and $x_{i\join j}$ pairwise distinct. (The case where exactly two of these coefficients are equal is analyzed below.)

If $\be_{-} = -1$, so that $2(\be_{-} + 1) =0$, and $\ep = 1$, then the column space of $((2\be_{-} - 1)\imat_{3} + \ones_{3}) = -3\imat_{3} + \ones_{3}$ is transverse to the right-hand side of \eqref{blocktalgidemdiffequationssimp2} and so the system \eqref{blocktalgidemdiffequationssimp2} is incompatible. Hence, if $\be_{-} \neq 1/2$ and $\ep = 1$, then the only remaining possibility for solutions of \eqref{blocktalgidemdiffequationssimp} is the case where exactly two of $x_{i}$, $x_{j}$, and $x_{i\join j}$ are equal. This is treated below, after consideration of the case $\be_{-} = 1/2$. 

If $\be_{-} = -1$ and $\ep = 0$, then a solution of \eqref{blocktalgidemdiffequationssimp2} has all its components equal, and this case has already been excluded.

If $\be_{-} = 1/2$, so that $\be = -\tfrac{n}{2(n-3)}$, then any $x_{i}$, $x_{j}$, and $x_{i\join j}$ such that $x_{i} + x_{j} + x_{i\join j} = \ep$ solve \eqref{blocktalgidemdiffequationssimp2}. Squaring this gives 
\begin{align}
\ep^{2} = x_{i}^{2} + x_{j}^{2} + x_{i\join j}^{2} + 2(x_{i}x_{j} + x_{i}x_{i\join j} + x_{j}x_{i\join j} ),
\end{align}
and substituting this and $\be_{+} + \be_{-} = 1 + 2\be$ into \eqref{blocktalgidemequationssimp} yields
\begin{align}\label{bc1}
\begin{split}
c &= x_{i}(x_{i} + x_{j} + x_{i\join j} - \ep) + 2\be(x_{i}x_{j} + x_{i}x_{i\join j} + x_{j}x_{i\join j} ) \\
&= 2\be(x_{i}x_{j} + x_{i}x_{i\join j} + x_{j}x_{i\join j} ) = \be(\ep^{2} -  x_{i}^{2} - x_{j}^{2} - x_{i\join j}^{2}).
\end{split}
\end{align}
On the other hand, by \eqref{nc},
\begin{align}\label{nc2}
\begin{split}
nc & = \tfrac{n(\be - 1)}{n-2}\left(\ep^{2} -  x_{i}^{2} - x_{j}^{2} - x_{i\join j}^{2} \right).
\end{split}
\end{align}
Comparing \eqref{bc1} and \eqref{nc2} shows 
\begin{align}
0 = 2((n-3)\be + 1)\left(\ep^{2} -  x_{i}^{2} - x_{j}^{2} - x_{i\join j}^{2} \right)=(n-2)\left(\ep^{2} -  x_{i}^{2} - x_{j}^{2} - x_{i\join j}^{2} \right),
\end{align}
so that 
\begin{align}
&x_{i} + x_{j} + x_{i\join j} = \ep, & &x_{i}^{2} + x_{j}^{2} + x_{i\join j}^{2} = \ep^{2}.
\end{align}
This shows \eqref{elel}.

There remains to consider the case in which $x_{i}$, $ x_{j}$, and $x_{i\join j}$ take exactly two distinct values. Write $x_{i} = x_{j} = s$ and $x_{i \join j} = t \neq s$. In this case, \eqref{blocktalgidemdiffequationssimp} becomes
\begin{align}\label{blocktalgidemdiffequationssimpsst}
\begin{split}
0 = &\resizebox{.85\textwidth}{!}{$\begin{pmatrix}
s-t & 0 & 0\\ 
0 & t-s & 0 \\ 
0 & 0 & 0
\end{pmatrix}
\left(
\left((2\be_{-} - 1)\imat_{3} + \ones_{3}\right)
\begin{pmatrix} s\\s\\t\end{pmatrix}
- \begin{pmatrix}\ep\\\ep\\\ep\end{pmatrix}
\right)
=\begin{pmatrix}
s-t & 0 & 0\\ 
0 & t-s & 0 \\ 
0 & 0 & 0
\end{pmatrix}
\begin{pmatrix} (2\be_{-} + 1)s + t - \ep\\ (2\be_{-} + 1)s + t- \ep\\2s +  (2\be_{-} + 1)t - \ep\end{pmatrix}$},
\end{split}
\end{align}
which forces
\begin{align}\label{sst1}
t = \ep - (2\be_{-} + 1)s.
\end{align}
The equation \eqref{blocktalgidemequationssimp} yields
\begin{align}\label{sst2}
c & = \left((1 + \be_{+} + \be_{-})s - \ep + 2\be_{+}t\right)s,\\
\label{sst3}c & = t^{2} - \ep t+ 2(\be_{+} + \be_{-}) st + (\be_{+} - \be_{-}) s^{2},
\end{align}
and \eqref{nc} becomes
\begin{align}\label{sst4}
nc & =(2 + 3\be_{+} + \be_{-})s^{2} - 2\ep s + 2(3\be_{+} + \be_{-})st + t^{2} - \ep t .
\end{align}
Substituting \eqref{sst1} in \eqref{sst2} yields
\begin{align}\label{sst5}
c & = s\left((1 - \be_{+} + \be_{-}  - 4\be_{+}\be_{-})s + \ep(2\be_{+} - 1) \right).
\end{align}
Subtracting \eqref{sst3} from \eqref{sst4}, substituting \eqref{sst1} in the result, and comparing with \eqref{sst5} yields
\begin{align}\label{sst6}
\begin{split}
(n-1)c & = 2(1 + \be_{+} + \be_{-})s^{2}  - 2\ep s + 4\be_{+}st\\
& = 2s\left((1 - \be_{+} + \be_{-} - 4\be_{+}\be_{-})s  + \ep(2\be_{+} - 1)\right) = 2c. 
\end{split}
\end{align}
so that $c = 0$, since it is assumed that $n > 3$. If $1 - \be_{+} + \be_{-} - 4\be_{+}\be_{-} = 0$, then $2\be_{+} \neq 1$, so that \eqref{sst6} implies $\ep s = 0$. If $\ep = 1$, then $s = 0$ and no new solution is obtained. If $\ep = 0$, then, by \eqref{sst1}, $t = -(2\be_{+} + 1)s$. In \eqref{sst3} this yields $0 = (1 + 3\be_{+} - 3\be_{-} - 4\be_{+}\be_{-})s^{2} = 4(\be_{+} - \be_{-})s^{2}$, so that $\be_{+} = \be_{-}$, which is inconsistent with $1 - \be_{+} + \be_{-} - 4\be_{+}\be_{-} = 0$. Hence there is no solution satisfying this last condition.
Suppose $1 - \be_{+} + \be_{-} - 4\be_{+}\be_{-} \neq 0$. 
In conjunction with \eqref{sst1} and \eqref{sst5}, \eqref{sst6} implies $s = 0$ and $t = \ep$ or 
\begin{align}\label{sst7}
\begin{split}
s & = \tfrac{\ep(1 - 2\be_{+})}{1 - \be_{+} + \be_{-} - 4\be_{+}\be_{-}}= \tfrac{\ep(n-2)(n-2(n-1)\be)}{ 4(n-3)(n-1)\be^{2} - 2n(n-4)\be + n(n-4)},\\
t & = \tfrac{\ep(\be_{+} - \be_{-})}{1 - \be_{+} + \be_{-} - 4\be_{+}\be_{-}} = \tfrac{2\ep(n-2)^{2}\be}{4(n-3)(n-1)\be^{2} - 2n(n-4)\be + n(n-4)}.
\end{split}
\end{align}
When $\ep = 1$, it follows that the three elements $e_{B}^{i}$, $e_{B}^{j}$ and $e_{B}^{i\join j}$ defined by \eqref{fbidem} are idempotents. When $\ep = 0$, no solutions are obtained.
\end{proof}

\begin{remark}
Because a Steiner triple system of order $n > 3$ has order $n \geq 7$ and, in this case,
\begin{align}
(n-2)^{2}( 1 - \be_{+} + \be_{-} - 4\be_{+}\be_{-}) = 4(n-3)(n-1)\be^{2} - 2n(n-4)\be + n(n-4)
\end{align}
is never zero in a Euclidean field, in this case \eqref{sst7} always has sense. 
\end{remark}

\begin{example}
Suppose $\B \in STS(n)$. When $\be = 1$, $e_{i}\mlt e_{j} = e_{i \join j}$. By Lemma \ref{talgidempotentlemma}, the elements $e^{0}_{B} = \tfrac{1}{3}(e_{i} + e_{j} + e_{i\join j})$, $e^{i}_{B} = \tfrac{1}{3}(2e_{i} - e_{j} - e_{i\join j})$,  $e^{j}_{B} = \tfrac{1}{3}(2e_{j} - e_{i} - e_{i\join j})$, and $e^{i \join j}_{B} = \tfrac{1}{3}(2e_{i\join j} - e_{i} - e_{j})$ associated with the block $B = \{i, j, i\join j\}$ are idempotents, and straightforward computations show they generate a subalgebra isomorphic to $\ealg^{3}(\fie)$.
\end{example}

\begin{example}
Suppose $\B \in STS(n)$. When $\be_{+} = 1/2$ (so $\be =\tfrac{n}{2(n-1)}$), for $B = \{i, j, i \join j\} \in \B$, the idempotents $e^{0}_{B}$, $e_{i}$, $e_{j}$, and $e_{i\join j}$ span a $3$-dimensional subalgebra. The nontrivial relations are $e^{0}_{B}\mlt e_{i} = \tfrac{1}{2}e^{0}_{B} + \tfrac{n-3}{2(2n-3)}e_{i}$ and those obtained from it with $j$ or $i\join j$ in place of $i$.
\end{example}

\begin{example}\label{notsimpleexample}
This example shows that $\talg_{1}(\fie, \B)$ need not be simple. Let $\B = AG(2, 3)$ be the affine plane of order $3$ as in \eqref{9124}. It is claimed that the algebra $\talg_{1}(\fie, \B)$ is isomorphic to the fourfold direct sum $\ealg^{2}(\fie)\oplus \ealg^{2}(\fie)\oplus \ealg^{2}(\fie)\oplus \ealg^{2}(\fie)$. (This observation was communicated to the author by Vladimir Tkachev.) 

Suppose $\be \in \fie \setminus \{-1/6\}$. 
Observe that $\B$ contains four maximal sets $\B_{I}$, $1 \leq I \leq 4$ of pairwise nonintersecting blocks, each having cardinality $3$ (for example $\{147, 258, 369\}$). They are given by the blocks corresponding respectively to rows, columns, left diagonals, and right diagonals. Because $\sum_{B \in \B_{I}}e^{0}_{B} = 0$, for distinct $B, B^{\prime} \in \B_{I}$, $e^{0}_{B}\mlt e^{0}_{B^{\prime}} = -e^{0}_{B} - e^{0}_{B^{\prime}}$. (As $\Aut(\B)$ acts transitively on the points and blocks of $\B$, to check a claim for all points and blocks it suffices to check it for any particular choice of points and blocks; for example for $B = 147$ and $B^{\prime} = 258$.) This suffices to show that $\balg^{e_{i}}_{i} = \fie\{e^{0}_{B}: B \in \B_{I}\}$ is a subalgebra of $\talg_{\be}(\fie, \B)$ isomorphic to $\ealg^{2}(\fie)$. 

That $B \in \B_{I}$ and $B^{\prime} \in\B_{J}$ with $i \neq j$ holds if and only if $B$ and $B^{\prime}$ are incident. There are four blocks $A$, $B$, $C$, and $D$ containing a given point $i \in \bar{9}$ (for example $\{123, 147, 168, 159\}$ are the blocks containing $1$) and a straightforward computation (for example with $A = 123$, $B = 147$, $C = 168$, $D = 159$) using $\ga_{A} + \ga_{B} + \ga_{C} + \ga_{D} = 3e_{i}$ shows
\begin{align}\label{affineplane1}
e^{0}_{A}\mlt e^{0}_{B} = \tfrac{1-\be}{6\be + 1}(e^{0}_{C} + e^{0}_{D}).
\end{align}
When $\be = 1$, \eqref{affineplane1} shows $e^{0}_{B} \mlt e^{0}_{B^{\prime}} = 0$, and this suffices to show $\balg^{e_{i}}_{i}\mlt \balg^{e_{i}}_{j} = \{0\}$, so that $\talg_{1}(\fie, \B) = \oplus_{i = 1}^{4}\balg^{e_{i}}_{i}$ is an orthogonal direct sum of subalgebras isomorphic to $\ealg^{2}(\fie)$, one for each maximal set of nonintersecting blocks.
\end{example}

\begin{example}\label{gabexample}
Consider a Hall triple system $\B \in HTS(n)$. For a block $B = \{i, j, i\join j\} \in \B$ and $k \notin B$, the set $k \join B$, defined by $k \join B = \{k\join i, k \join j, k \join (i\join j)\}$, is again a block of $\B$ by \eqref{hallid1}. 

Let $B = \{i, j, i\join j\} \in \B$ and $k \notin B$. Note that $i\join (k \join B) = (i\join k)\join B = (i \join j)\join (k \join B) = (j\join k)\join B = j\join (k \join B)$. Using \eqref{hallid1} and \eqref{hallid2} it is straightforward to check the following relations:
\begin{align}
\begin{aligned}
\ga_{B}&\mlt e_{i} = (1 + \be_{-})e_{i} + \be_{+}\ga_{B},\\
\ga_{B}&\mlt (e_{i} - e_{j}) = (1 + \be_{-})(e_{i} - e_{j}),\\
\ga_{B}&\mlt e_{k} = \al \ga_{B} + 3\al e_{k} + \be \ga_{k \join B},\\
\ga_{B}&\mlt(e_{i\join k} - e_{j\join k}) = \tfrac{3(\be - 1)}{n-2}(e_{i\join k} - e_{j\join k}),\\
\ga_{B}&\mlt \ga_{k\join B}  = 3\al(\ga_{B} + \ga_{k \join B}) + 3\be \ga_{(i\join k)\join B},\\
\ga_{B}&\mlt(\ga_{k\join B} - \ga_{(i \join k) \join B)}) = 3\be_{-}(\ga_{k\join B} - \ga_{(i \join k) \join B)}).
\end{aligned}
\end{align}
This shows that the eigenvalues of $L(\ga_{B})$ include $1 + 3\be_{+} + \be_{-} = \tfrac{2n\be + n-6}{2n}$, with multiplicity at least $1$; $1+ \be_{-} = \tfrac{(n-3)(1-\be)}{n-2}$, with eigenspace containing the two-dimensional subspace $\{\la_{i}e_{i} + \la_{j}e_{j} + \la_{i\join j}e_{i\join j}: \la_{i} + \la_{j} + \la_{i\join j} = 0\}$; $\tfrac{3(\be - 1)}{n-2}$, with eigenspace containing the subspace $\fie\{e_{k\join i} - e_{k\join j}: k \notin B\}$; and $3\be_{-} = -\tfrac{3((n-3)\be + 1)}{n-2}$, with eigenspace containing the nontrivial subspace $\fie\{\ga_{k \join B} - \ga_{(i \join k)\join B}, k \notin B\}$. For $n > 9$ it is not clear that this list is exhaustive, but when $n = 9$, because the automorphism group of $AG(2, 3)$ (up to isomorphism the unique element of $STS(9)$) is transitive on points and blocks, it suffices to check that the eigenvalues and eigenspaces of $L(e^{0}_{B})$ are as indicated in Table \ref{gabtable} (supposing $\be \neq -1/6$).
\begin{table}
\begin{tabular}{|c|c|c|}
\hline
Eigenvalue &  Elements spanning eigenspace & Multiplicity\\
\hline
$1$ & $e^{0}_{123}$ & $1$ \\
\hline
$\tfrac{2(1-\be)}{6\be + 1}$ & $e_{1} - e_{2}$, $e_{2} - e_{3}$ & $2$\\
\hline
$\tfrac{\be-1}{6\be +1}$ & 
\makecell{
$e_{4\join 1} - e_{4\join 2} = e_{7} - e_{9}$, 
$e_{5\join 1} - e_{5\join 2} = e_{9} - e_{8}$, \\
$e_{7\join 1} - e_{7\join 2} = e_{4} - e_{6}$, 
$e_{8\join 1} - e_{8\join 2} = e_{6} - e_{5}$, }
 & $4$\\
\hline
$-1$ & $e^{0}_{4 \join 123} - e^{0}_{7\join 123} = e^{0}_{789} - e^{0}_{456}$ & $1$\\
\hline
\end{tabular}
\caption{Eigenvalues of $L(e^{0}_{123})$ in $(\talg_{\be}(\fie, AG(2, 3)), \mlt)$.}\label{gabtable}
\end{table}
The indicated multiplicities are incorrect in certain special cases where the indicated eigenvalues coincide.
\end{example}

\begin{lemma}\label{aglemma}
For $\be \in \fie \setminus\{1, -1/6\}$ such that $6\be^{2} +1 \neq 0$, the algebra $(\talg_{\be}(\fie, AG(2, 3)), \mlt)$ is simple.
\end{lemma}
\begin{proof}
The content of the lemma is that $(\talg_{\be}(\fie, AG(2, 3)), \mlt)$ is simple for the value $\be = -4/3$ excluded in the statement of Theorem \ref{stsalgebrasimpletheorem}. Write $\B = AG(2, 3)$. Because $\sum_{i \in B \in \B}\ga_{B} = re_{i} = \tfrac{n-1}{2}e_{i}$, the set $\axes = \{e^{0}_{B}: B \in \B\}$ comprises idempotent spanning $\talg_{\be}(\fie, \B)$.
By Example \ref{gabexample}, for any $B \in \B$, the endomorphism $L(e_{B}^{0})$ is diagonalizable with eigenvalues $1$, $-2/3$, $1/3$, and $-1$, having multiplicities $1$, $2$, $4$, and $1$. Because $6\be^{2} +1 \neq 0$, by Theorem \ref{stsalgebratheorem}, $\kform$ is nondegenerate and invariant, the set $\axes$ and the form $\kform$ satisfy the hypotheses of Lemma \ref{simplelemma}, and so a nontrivial ideal $\ideal \subset \talg_{\be}(\fie, \B)$ contains $e^{0}_{B}$ for some $B \in \B$. Because $L(e_{B}^{0})$ is invertible, $\ideal = \talg_{\be}(\fie, \B)$.
\end{proof}

The same argument would prove the simplicity of $(\talg_{\be}(\fie, \B), \mlt)$ for $\B \in HTS(n)$ and $\be = -\tfrac{n-1}{n-3}$ were it possible to establish that the multiplicity of $1$ as an eigenvalue of $L(e_{B}^{0})$ is $1$ and $L(e_{B}^{0})$ is invertible. When $n > 9$ this appears to require more serious consideration of the combinatorial structure of $\B$.

It would be interesting to resolve fully the following two questions:
\begin{itemize}
\item For $\B \in STS(n)$ and $\be \neq \bar{\be} \in \fie$ decide if $\talg_{\be}(\fie, \B)$ and $\talg_{\bar{\be}}(\fie, \B)$ are isomorphic. 
\item For $\be \in \fie$ and distinct $\B, \B^{\prime} \in STS(n)$ decide if $\talg_{\be}(\fie, \B)$ and $\talg_{\be}(\fie, \B^{\prime})$ are isomorphic.
\end{itemize}
This appears to require more detailed information about idempotents answering questions such as: if an idempotent $e \in \talg_{\be}(\fie, \B)$ satisfies $\kform(e, e) = \kform(e_{i}, e_{i})$ for all $i \in \setn$, is there $j \in \setn$ such that $e = e_{j}$? Lemma \ref{permutelemma} answers this question affirmatively, but under quite restrictive hypotheses both on $\fie$ and the values of $\be$. It is not clear how to answer even something so apparently simples as: is it the case that $\talg_{\be}(\fie, \B)$ contains no nontrivial square-zero element if $\be \neq \tfrac{6-n}{2n}$?

\subsection*{Acknowledgments}
I thank Vladimir Tkachev for discussions related to the contents of this paper and for sharing his related work in progress.

\bibliographystyle{amsplain}
\def\polhk#1{\setbox0=\hbox{#1}{\ooalign{\hidewidth
  \lower1.5ex\hbox{`}\hidewidth\crcr\unhbox0}}} \def\cprime{$'$}
  \def\cprime{$'$} \def\cprime{$'$}
  \def\polhk#1{\setbox0=\hbox{#1}{\ooalign{\hidewidth
  \lower1.5ex\hbox{`}\hidewidth\crcr\unhbox0}}} \def\cprime{$'$}
  \def\cprime{$'$} \def\cprime{$'$} \def\cprime{$'$}
  \def\polhk#1{\setbox0=\hbox{#1}{\ooalign{\hidewidth
  \lower1.5ex\hbox{`}\hidewidth\crcr\unhbox0}}} \def\cprime{$'$}
  \def\Dbar{\leavevmode\lower.6ex\hbox to 0pt{\hskip-.23ex \accent"16\hss}D}
  \def\cprime{$'$} \def\cprime{$'$} \def\cprime{$'$} \def\cprime{$'$}
  \def\cprime{$'$} \def\cprime{$'$} \def\cprime{$'$} \def\cprime{$'$}
  \def\cprime{$'$} \def\cprime{$'$} \def\cprime{$'$} \def\dbar{\leavevmode\hbox
  to 0pt{\hskip.2ex \accent"16\hss}d} \def\cprime{$'$} \def\cprime{$'$}
  \def\cprime{$'$} \def\cprime{$'$} \def\cprime{$'$} \def\cprime{$'$}
  \def\cprime{$'$} \def\cprime{$'$} \def\cprime{$'$} \def\cprime{$'$}
  \def\cprime{$'$} \def\cprime{$'$} \def\cprime{$'$} \def\cprime{$'$}
  \def\cprime{$'$} \def\cprime{$'$} \def\cprime{$'$} \def\cprime{$'$}
  \def\cprime{$'$} \def\cprime{$'$} \def\cprime{$'$} \def\cprime{$'$}
  \def\cprime{$'$} \def\cprime{$'$} \def\cprime{$'$} \def\cprime{$'$}
  \def\cprime{$'$} \def\cprime{$'$} \def\cprime{$'$} \def\cprime{$'$}
  \def\cprime{$'$} \def\cprime{$'$} \def\cprime{$'$}
\providecommand{\bysame}{\leavevmode\hbox to3em{\hrulefill}\thinspace}
\providecommand{\MR}{\relax\ifhmode\unskip\space\fi MR }
\providecommand{\MRhref}[2]{%
  \href{http://www.ams.org/mathscinet-getitem?mr=#1}{#2}
}
\providecommand{\href}[2]{#2}

\end{document}